\documentclass{cedram-aif}

\usepackage{amsmath,amssymb,amsthm,amscd}
\usepackage{graphicx}
\usepackage{mathrsfs}
\usepackage[mathcal]{eucal}

\def\arxivno#1{{\texttt{#1}}}

\numberwithin{equation}{section}

\newtheorem{Thm}{Theorem}%[section]
\newtheorem{Cor}{Corollary}

\newtheorem{Lem}{Lemma}

\theoremstyle{definition}
\newtheorem{Def}{Definition}

\theoremstyle{remark}
\newtheorem{Rem}{Remark}
\newtheorem{Ex}{Example}

\edef\Loja.{\L ojasiewicz}

\renewcommand\ge{\geqslant}
\renewcommand\le{\leqslant}
\let\tildeaccent=\~
\let\hataccent=\^
\renewcommand\~[1]{\widetilde{#1}}
\renewcommand\^[1]{\widehat{#1}}

\def\<{\left<}
\def\>{\right>}
\def\({\left(}
\def\){\right)}

\def\norm#1{\left\Vert #1 \right\Vert}

\let\parasymbol=\S
\def\secref#1{\parasymbol\ref{#1}}
\let\ssm=\smallsetminus

\let\simeq=\cong

\def\Var{\operatorname{Var}}

\def\Mat{\operatorname{Mat}}

\def\Re{\operatorname{Re}}
\def\Im{\operatorname{Im}}
\def\Arg{\operatorname{Arg}}
\def\dist{\operatorname{dist}}

\def\etc.{\emph{etc}.}

\def\:{\colon}
\def\R{{\mathbb R}}
\def\C{{\mathbb C}}
\def\Z{{\mathbb Z}}
\def\N{{\mathbb N}}

\def\e{\varepsilon}
\def\S{\varSigma}
\def\l{\lambda}

\def\pf{\hbox{$P$\kern-6pt\raise-2pt\hbox{$F$}}}
\def\f{\varphi}

\def\res{\operatornamewithlimits{Res}}
\def\d{\,\mathrm d}
\def\iu{{\mathrm i}}

\def\cN{\mathcal{N}}

\def\cI{\mathcal{I}}
\def\clF{\mathscr F_{n,m}}
\def\clS{\mathscr S_{n,m}}
\def\slope{\operatorname{\angle}}
\def\Q{{\mathbb Q}}

\begin{document}

\title[Polynomial bounds for oscillation of Fuchsian systems]
 { Polynomial bounds for oscillation \\of solutions of Fuchsian systems  }

 \dedicatory{To Ma\hataccent\i tre Bernard for his eightieth birthday}

\author
 {\firstname{Gal} \lastname{Binyamini}
 \and
 \firstname{Sergei} \lastname{Yakovenko}
 }

\begin{abstract}
We study the problem of placing effective upper bounds for the
number of zeros of solutions of Fuchsian systems on the Riemann
sphere. The principal result is an explicit (non-uniform) upper
bound, polynomially growing on the frontier of the class of Fuchsian
systems of a given dimension $n$ having $m$ singular points. As a
function of $n,m$, this bound turns out to be  double exponential in
the precise sense explained in the paper.

As a corollary, we obtain a solution of the so called restricted
infinitesimal Hilbert 16th problem, an explicit upper bound for the
number of isolated zeros of Abelian integrals which is polynomially
growing as the Hamiltonian tends to the degeneracy locus. This
improves the exponential bounds recently established by A.~Glutsyuk
and Yu.~Ilyashenko.
\end{abstract}

\address{Weizmann Institute of Science\\Rehovot 76100\\ISRAEL}
\email{\{gal.binyamini, sergei.yakovenko\}@weizmann.ac.il}


\maketitle

\section{Zeros of solutions of Fuchsian systems and restricted infinitesimal Hilbert
16th problem}

\subsection{Fuchsian systems and zeros of their solutions} Let $\Omega$ be
a meromorphic (rational) $n\times n$-matrix 1-form on the Riemann sphere
$\C P^1$ with a singular (polar) locus $\S=\{\tau_1,\dots,\tau_m\}$ consisting
of $m$ distinct points. The linear system of Pfaffian equations
\begin{equation}\label{eq:ls}
    \d x-\Omega x=0,\qquad
    \Omega=
    \begin{pmatrix}
    \omega_{11}&\cdots&\omega_{1n}\\\vdots&\ddots&\vdots\\\omega_{n1}&\cdots&\omega_{nn}
    \end{pmatrix},
    \quad
    x=
    \begin{pmatrix}
    x_1\\\vdots\\x_{n}
    \end{pmatrix},
\end{equation}
is said to be \emph{Fuchsian}, if all points $\tau_i\in\S$ are
first order poles of $\Omega$. We will denote the collection of
all Fuchsian systems of rank $n$ with $m$ singular points by
$\clF $. Later it will be identified with an
$m(n^2+1)$-dimensional semialgebraic variety. A Fuchsian system
can be always expanded as the sum of its \emph{principal
Fuchsian parts},
\begin{equation}\label{fuchs-principal-parts}
    \Omega=\sum_{j=1}^m\Omega_j,\qquad
    \Omega_j=\frac{A_j\d t}{t-\tau_j},\quad \tau_j\in\S,
    \quad \sum_{j=1}^m A_j=0,
\end{equation}
with the \emph{residue matrices} $A_1,\dots,A_m$ (the point
$\tau=\infty$, singular for each $\Omega_j$, is nonsingular for the
sum; if $\infty=\tau_m\in\S$, then one has $m-1$ terms corresponding
to the finite points $\tau_1,\dots,\tau_{m-1}$ of $\S$, yet with a
nontrivial residue $A_m$ corresponding to the point at infinity).

Solutions of a Fuchsian system~\eqref{eq:ls} are multivalued
holomorphic vector-functions on $\C P^1\ssm\S$, ramified over the
singular locus, and growing moderately (no faster than polynomially
in the reciprocal distance to $\S$) when approaching points of $\S$
along non-spiraling paths.

We are interested in \emph{explicit global upper bounds} on the
maximal number of isolated zeros of linear combinations of the form
$y(t)=\sum_{i=1}^n c_i\,x_i(t)$, $t\in\C P^1$. The initial
motivation for this question came from the attempts to solve the so
called infinitesimal Hilbert 16th problem, see
\cite{montreal,centennial} and \secref{sec:inf-16}, but the subject
is interesting in itself as well, as linear systems provide a vast and
essentially unique source of transcendental functions which still
admit global finite bounds for their zeros \cite[\parasymbol 3]{montreal}.

Because of the multivaluedness of the solutions, some disambiguation
is required. Assume that the punctured Riemann sphere $\C P^1\ssm\S$
is covered by finitely many \emph{simply connected polygonal
domains} $U_\alpha$ (polygonality here and below means that each
domain $U_\alpha$ is bounded by finitely many circular arcs and line
segments). In each such domain one can consistently choose a
holomorphic branch of any solution $x(t)$ to \eqref{eq:ls} and form
all scalar linear combinations $y(t)$ as above. Clearly, the number
     $$N_\alpha(\Omega)=\sup_{x(\cdot),\ c\in\C^n}
      \#\{t\in {U_\alpha}\: c_1 x_1(t)+\cdots+c_n
      x_n(t)=0\}\le+\infty,$$
where $\#\{\text{a set}\}$ denotes the number of \emph{isolated}
points of the set and the supremum is taken over all solutions
$x(\cdot)$ of the system \eqref{eq:ls}, is independent of the choice
of the branches.

Let $N(\Omega,U)$ be the total number of the zeros counted in this
way in all domains of a given covering $U=\{U_\alpha\}$, and
$\cN(\Omega)$ the infimum of this total number, taken over all
polygonal coverings $U$,
\begin{equation} \label{eq:counting-function}
\begin{gathered}
    N(\Omega,U)=\sum_\alpha N_\alpha(\Omega)\le+\infty,
    \\
    \cN(\Omega)=\inf_{U}N(\Omega,U).
    \end{gathered}
\end{equation}
%To avoid trivial counterexamples, we need to assume that the
%total number of arcs bounding all the domains $U_\alpha$, is
%commensurable with the homology of $U$ (measured as the number
%of singular points of $\Omega$). For our purposes it will be
%sufficient to assume that the number of arcs grows no faster
%than linearly with $m$. It is clear that this condition is of a
%purely technical nature: we could alternatively require that
%$U$ is a triangulation of $\C P^1\ssm\S$ into a minimal number
%of open triangles bounded by real segments. The bounds, while
%literally different for two different constructions, will be
%obviously equivalent in the sense of explicit asymptotic growth
%(see~\secref{sec:semialgebraic-bounds} below). \NB{Why is this
%paragraph needed, when we take \emph{infimum}?} \NB{Since we
%theoretically could decrease the number of zeros by considering
%polygonal approximations to spiralling domains...}

By abuse of language, we say that the number $\cN(\Omega)$, if
finite, gives a \emph{global bound for the number of zeros of
solutions of the Fuchsian system} \eqref{eq:ls}. Considered as a
function  $\cN\:\clF \to\N\cup\{+\infty\}$ on the space of all
Fuchsian systems (or its subspaces), it will be called the
\emph{counting function}.

Without additional assumptions the ``bound'' $\cN(\Omega)$ may
well be infinite, since even individual linear combinations may
have infinitely many isolated zeros accumulating to the
singular locus $\S$. It turns out that the finiteness of the bound is
closely related to the spectral properties of the \emph{residue
matrices}
  $$
  A_i=\res_{\tau_i}\Omega\in\Mat_{n\times n}(\C), \qquad \tau_i\in\S,
  \quad i=1,\dots,m.
  $$

\begin{Ex}\label{ex:nonreal-spec}
The function $y(t)=t^\iu+t^{-\iu}=2\cos\ln t$ is a linear
combination corresponding to the \emph{Euler system} with
$\Omega=t^{-1}\bigl(\begin{smallmatrix}\phantom{-}0&1\\-1&0\end{smallmatrix}\bigr)
\d t$ in the standard affine chart on $\C P^1$. This combination has
infinitely many zeros accumulating to the two singular points
$\tau_1=0$ and $\tau_2=\infty$ with the residues
$A_0=\bigl(\begin{smallmatrix}\phantom{-}0&1\\-1&0\end{smallmatrix}\bigr)=-A_\infty$.
\end{Ex}

However, under additional assumptions on the spectra of the
residue matrices $A_i$ one can guarantee finiteness of the
upper bound. The following result can easily be obtained by the
tools developed in \cite{asik:finiteness,asik:fewnomials}.

\begin{Thm}\label{thm:quasivk}
Assume that all residue matrices $A_i=\res_{\tau_i}\Omega$ of the Fuchsian
system \eqref{eq:ls} have only real eigenvalues,
\begin{equation*}
    \operatorname{Spec}A_i\subseteq\R,\qquad i=1,\dots,m.
\end{equation*}
Then the corresponding bound for the number of isolated zeros of solutions
is finite, $\cN(\Omega)<\infty$.
\end{Thm}

We will not give the proof of this result since it is purely
existential and does not allow any explicit estimate for the bound
$\cN(\Omega)$. Formally it follows from Theorem~\ref{thm:main}
below, whose proof is independent of Theorem~\ref{thm:quasivk}.

\subsection{Unboundedness of the counting function}
Denote by $\clS \subset\clF $ the subspace of the space of Fuchsian
systems with real spectra of all residues:
\begin{equation}\label{eq:spec}
    \clS =\{\Omega\in\clF \:\forall
    \tau\in\S,\
    \operatorname{Spec}\res_\tau\Omega\subseteq\R\}.
\end{equation}
By Theorem~\ref{thm:quasivk}, the counting function $\cN$ takes
finite values on $\clS $. Yet since $\clS$ is non-compact, this
finiteness does not prevent the counting function from being
\emph{unbounded} on $\clS $. The following examples show that this
is indeed the case.

\begin{Ex}\label{ex:norm-dependence}
The counting function may grow to infinity together with the norms
of residue matrices of $\Omega$.

Indeed, for the Euler system,
$\Omega=t^{-1}\bigl(\begin{smallmatrix}N&0\\0&0\end{smallmatrix}\bigr)\d{}t$
with $N\in\N$, the linear combination $y(t)=t^N-1$ has $N$ isolated zeros uniformly
distributed on the unit circle. In this example $N$ is the norm of
the residues at $t=0$ and $t=\infty$.
\end{Ex}

Somewhat less expected is the fact that the counting function may
grow to infinity as some of the singular points collide with each
other even without explosion (growth to infinity) of the residues.

\begin{Ex}\label{ex:distance-dependence}
Consider the family of Fuchsian systems $\Omega_\e$, $\e\in\C$,
$\e\ne0$, with the four-point singular locus
$\S=\{0,\e,1/\e,+\infty\}$ and four residue matrices, each with the
real spectrum,
\begin{equation}\label{eq:nonreal-limit-spec}
    A_0=\begin{pmatrix} 0&1\\ & 0\end{pmatrix},
    \
    A_\e=\begin{pmatrix} \phantom{-}0&\\ -1 &0 \end{pmatrix},
    \qquad
    \begin{aligned}
    A_{1/\e}&=-A_\e,\\ A_\infty&=-A_0.
    \end{aligned}
\end{equation}
As the complex parameter $\e\ne0$ tends to zero, two pairs of
singularities, $\{0,\e\}$ and $\{1/\e,\infty\}$, collide. The
limit is the Euler system described in
Example~\ref{ex:nonreal-spec} with the residue matrices
$A_0+A_\e$ and $-(A_\infty+A_{1/\e})$ respectively, which has
infinitely many isolated zeros accumulating to both
singularities.

This obviously means that
$\limsup_{|\e|\to0}\cN(\Omega_\e)=+\infty$. Indeed, otherwise there
would exist a finite upper bound $N$ for $\cN(\Omega_\e)$ uniform
over all sufficiently small $|\e|$. Consider a triangle $K\Subset\C
P\ssm\S$ which contains at least $N+1$ isolated zeros of the
solution $y(t)=t^\iu+t^{-\iu}$ of the limit system $\Omega_0$
\emph{strictly inside} (cf.~with Example~\ref{ex:nonreal-spec}). By
semicontinuity, for all sufficiently small values of $|\e|$, the
corresponding combination of solutions of the system with the matrix
$\Omega_\e$ will have at least $N+1$ zeros, since the restriction of
$\Omega_\e$ on $K$ depends holomorphically on $\e$. Hence
$\cN(\Omega_\e)\ge N+1$ in contradiction with the assumption.
\end{Ex}

\subsection{Polynomial bounds}
The two examples above show that the counting function $\cN(\Omega)$
grows to infinity as $\Omega$ approaches the frontier of the
subclass $\clS $ of Fuchsian systems, at least along some parts of
this frontier. The main positive result of the current work is an
\emph{explicit upper bound} for $\cN(\Omega)$ which grows
\emph{polynomially} as the Fuchsian system~\eqref{eq:ls}, identified
with its matrix $\Omega$, approaches the frontier of the class $\clS
$.

Here and everywhere in this paper we define the norm of a
polynomial $p\in\C[x_1,\dots,x_n]$ in several variables as follows,
\begin{equation}\label{eq:poly-norm}
    \|p\|=\sum |c_\alpha|,\qquad\text{if}
    \quad p(x)=\sum c_\alpha x^\alpha,\quad \alpha=(\alpha_1,\dots,\alpha_r)\in\Z_+^r.
\end{equation}
The norm $|M|$ of a matrix $M\in\Mat_{n\times n}(\C)$ is taken to be
the standard Euclidean (Hermitian) matrix norm (though any other
norm would of course be essentially equivalent).

Let $R_\flat\:\clF \to\R_+$ be the function defined as follows,
\begin{equation}\label{eq:main-R}
   R_\flat(\Omega)=2+\sum_{\tau'\ne\tau''\in\S}{\dist^{-1}(\tau',\tau'')}+
     \sum_{\tau\in\S}|\res\nolimits_{\tau}\Omega\,|,
\end{equation}
where $\dist(\cdot,\cdot)$ is the Fubini--Study distance on $\C
P^1$. The function $R_\flat$ serves to measure the (reciprocal)
distance to the frontier of the set $\clF $. Here and below $\R_+$
is a closed subset of \emph{nonnegative} real numbers.

\begin{Thm}\label{thm:main}
The counting function $\cN(\cdot)$ for the number of isolated zeros
of solutions of Fuchsian systems is explicitly bounded on the subset
$\clS \subset\clF $ by a power of the function $R_\flat(\cdot)$ as
follows\textup:
\begin{equation}\label{mainres}
    \cN\big|_{\clS }\le R_\flat^\nu,\qquad
      \nu=\nu_{n,m}\le {2^{O(n^2m)}}.
\end{equation}
The constant in $O(n^2m)$ is explicit and computable.
\end{Thm}

\begin{Rem}\label{rem:non-conf}
Note that the Fubini-Study distance on $\C P^1$ is not \emph{conformally}
invariant, whereas the counting function $\cN$ obviously is. This may be
used to improve the bounds: for instance, if $m=3$, then any three points
can be placed at a distance 1 from each other by a suitable conformal change
of variable $t$.
\end{Rem}

\begin{Rem}
The term $O(n^2m)$ appearing in the bound above is asymptotically
equivalent to the dimension of the variety $\clS$. A careful
inspection of the proof will reveal that this quantity in fact plays
the key role in determining the magnitude of the upper bound. On
suitably defined subsets one may apply the same techniques to
establish similar bounds depending on the dimension of the set,
though we don't explore this direction here.
\end{Rem}

The proof of Theorem~\ref{thm:main} is based on the controlled
process of reduction of the system \eqref{eq:ls} to a scalar $n$th
order differential equation along the lines suggested by
A.~Grigoriev \cite{alexg:thesis,alexg:arxiv}, see also
\cite{ongrig}. As a result, we obtain an explicit uniform upper
bound for the variation of argument of all linear combinations along
arcs sufficiently distant from the singular locus of this
equation.

By the argument principle, this implies an explicit upper bound
for the number of isolated zeros away from the singular locus
$\S$ of the initial system. As for the zeros arbitrarily close
to $\S$, we use the method developed in \cite{mrl-96} for
Fuchsian equations.

\begin{small}
\subsection{Applications for the infinitesimal Hilbert 16th
problem}\label{sec:inf-16} The principal motivation for
Theorem~\ref{thm:main} was the so called \emph{restricted
infinitesimal Hilbert sixteenth problem} as it was formulated
in~\cite{iglu:restricted,ilyash:openprob}.

Let $H\in\R[x,y]$ be a real polynomial of degree $n+1$ and
$\omega=p\,\d x+q\,\d y$ a real polynomial 1-form of degree
$n=\max(\deg p,\deg q)$. The problem is to find explicit upper
bound for the number $I_{\omega,H}$ of real ovals
$\delta_t\subset\{H=t\}$ on the \emph{nonsingular} affine level
curves of $H$, such that the Abelian integral
\begin{equation}\label{eq:abint}
    y_{\omega,H}(t)=\oint_{\delta_t}\omega, \qquad t\in (a,b),
\end{equation}
vanishes (assuming that such ovals are isolated). The supremum
\begin{equation}\label{eq:abint-bd}
    \mathcal I(H)=\sup_{\deg\omega\le n}
    I_{\omega,H}
\end{equation}
depends only on $H$ and is always finite, so that $\mathcal I$
can be considered as a function on the linear space $\mathscr
H_n=\{\deg H\le n+1\}$ of real polynomials of the specified
degree.

In \cite{asik:finiteness,varchenko:finiteness} it was shown that the
value $\cI(H)$ admits a uniform  bound,
\begin{equation}\label{eq:vk}
    \sup_{H\in\mathscr H_n}\cI(H)=\mathcal I_n<+\infty,
\end{equation}
yet the proof is non-constructive and the growth of $\mathcal
I_n$ as a function of $n$ is completely uncontrollable.

All known explicit bounds for $\cI(H)$ are non-uniform and
finite only for a residual (open dense) subset of polynomials
$H$ for $n\ge 3$. The only case where a uniform bound was known
to be extendable to a class of degenerate Hamiltonians, is that
of hyperelliptic integrals, see \cite{era-99}. Yet the bound in
this case is given by a tower function (iterated exponent) of
unspecified height (iteration depth).

As a corollary of Theorem~\ref{thm:main} we obtain a
\emph{polynomial} global upper bound for the number of zeros of
Abelian integrals on the same set $\mathscr M_n$. More precisely,
consider the expansion of $H$ as the sum of the principal
homogeneous parts,
$$
 H=H_{n+1}+\cdots+H_1+H_0,\qquad H_k(\l x,\l y)=\l^k H_k(x,y).
$$
Let $\mathscr M_n$ be a semialgebraic subset in $\mathscr H_n$,
defined by the following conditions,
\begin{enumerate}
 \item the principal homogeneous part $H_{n+1}$ is \emph{nonzero}:
 $$H_{n+1}(x,y)=\sum_{i+j=n+1} c_{ij} x^i y^{j}\not\equiv0;$$
 \item the principal part is square-free, i.e., $H_{n+1}$
     factors as a product of $n+1$ pairwise different
     linear terms corresponding to different points
     $\zeta_1,\dots,\allowbreak \zeta_{n+1}$ of the
     projective line $\C P^1$;
 \item the polynomial $H$ itself has exactly $n^2$ different critical points
 in $\C^2$ with pairwise different critical values $\tau_1,\dots,\tau_{n^2}\in\C$
 (the corresponding critical points will automatically be
 nondegenerate, i.e., the Hessian of $H$ is nonvanishing at each of them).
\end{enumerate}
Obviously, $\mathscr H_n$ is a complex affine space and $\mathscr
M_n$ a semialgebraic subset in this space. Let $R_\natural\:\mathscr
M_n\to\R_+$, $H\mapsto R_\natural(H)$ be the function defined in
terms of the homogeneous decomposition of $H$ as follows,
\begin{equation}\label{eq:carp-morse}
    R_\natural(H)=2+\frac{\|H\|}{\|H_{n+1}\|}+\sum_{i\ne j}^{n+1}\dist^{-1}(\zeta_i,\zeta_j)
    +\sum_{j=1}^{n^2}|\tau_j|+\sum_{i\ne j}^{n^2}|\tau_i-\tau_j|^{-1}
\end{equation}
This is a continuous semialgebraic function on the set  $\mathscr
M_n$ with the infinite limit on the frontier $\partial\mathscr M_n$.

\begin{Thm} \label{thm:abelian-main}
The number of isolated zeros of Abelian integrals admits a computable semialgebraic
bound on $\mathscr M_n$\textup:
\begin{equation}\label{bound-H}
    \cI(H)\big|_{\mathscr M_{n}}\le R_\natural^\nu,\qquad
    \nu=\nu_{n}\le {2^{O(n^4)}},
\end{equation}
As in Theorem~\ref{thm:main}, the constant in $O(n^4)$ is explicit and computable.
\end{Thm}

\begin{Rem}
Very recently A.~Glutsyuk and Yu.~Ilyashenko in a series of
papers
\cite{iglu:restricted,iglu:dan,glu:determinant,glu:topology}
achieved another bound for the number of zeros of Abelian
integrals, based on completely different ideas. Their bound is
finite \emph{on the same residual subset} $\mathscr M_n$, yet
for any fixed $n$ it grows \emph{exponentially} near some parts
of the boundary $\partial \mathscr M_n\subset\mathscr H_n$.
\end{Rem}

Theorem~\ref{thm:abelian-main} is a corollary to
Theorem~\ref{thm:main} and follows immediately from the fact that
Abelian integrals satisfy a hypergeometric Picard--Fuchs system of
linear ordinary differential equations explicitly constructed in
\cite{redundant}. The interested reader may easily restore the
details of the derivation.

We do not give a detailed proof of this Theorem here, since in a
separate publication \cite{inf16} we establish a \emph{uniform}
explicit bound for the number of zeros of Abelian integrals, thus
solving the \emph{unrestricted} infinitesimal Hilbert 16th problem.
\par

\end{small}

\subsection{Acknowledgements}
The authors are grateful S. Basu,  A. Gabrielov, A. Khovanskii
and N. Vorobjov. A special gratitude is due to Dmitry Novikov
for his invaluable assistance and endless discussions which
helped put the exposition in its current form. The research was
partially supported by the Minerva Foundation. One of the
authors (S.Y.) is incumbent of the Gershon Kekst professorial
chair.

\section{On the general nature of non-uniform bounds}
 \label{sec:semialgebraic-bounds}

The form of the bound established in Theorem~\ref{thm:main} begs for
further analysis, as one can question the form in which the
numerous different parameters determining the Fuchsian system
\eqref{eq:ls} are incorporated into a single function
$R_\flat(\Omega)$. In this section we discuss the general form of
explicit \emph{non-uniform} bounds on semialgebraic parameter
spaces, of which \eqref{mainres} is a particular case, and study
their universality.

This allows us to conclude that the double exponential
dependence of the bound \eqref{mainres} is independent of the
particular form chosen for $R_\flat$, as long as certain
``natural'' conditions are satisfied. Furthermore, we show that
the freedom of choosing the ``natural'' inverse distance
$R_\flat$ to the frontier will affect the powers $\nu_n$,
changing them by a polynomial in $n$ and $m$. Thus the
exponential bound for the powers $\nu_{n,m}$ (double
exponential with respect to $R_\flat$) \emph{cannot be
substantially improved in the class of non-uniform bounds on
the whole class of Fuchsian systems}. We explain the precise
meaning of the ``natural bounds'' in this section.

\subsection{Semialgebraic sets and their carpeting functions}
Recall that a subset $Z\subseteq\R^N$ of the real affine space is
called \emph{semialgebraic}, if it is a finite union of subsets each
defined by finitely many polynomial equalities and inequalities of
the form $\{f_\alpha(x)=0,\ g_\beta(x)<0\}$, where
$f_\alpha,g_\beta\in \R[x_1,\dots,x_N]$ are real polynomials. A
projective subset $Z\subseteq\R P^N$ is semialgebraic, if it is
semialgebraic in some (hence any) affine chart. A subset in $\C^N$ is
semialgebraic, if it is semialgebraic in the corresponding
``realification'' $\R^{2N}$; \etc.

The category of semialgebraic sets is stable by set theoretic
operations (finite union, intersection, complement) as well as by
projections (as asserted by the celebrated Tarski--Seidenberg
principle). As a result, we conclude that any formula involving
quantifiers and polynomial expressions, defines a semialgebraic
set. For the same reasons the image of a semialgebraic set by a
polynomial or semialgebraic map is again semialgebraic (a map
is semialgebraic if its graph is a semialgebraic subset of the
Cartesian product of the domain and the range).

The space $\clF$ of Fuchsian systems can be identified in several
ways with semialgebraic subsets in suitable affine spaces, see
\secref{sec:prod-struct} and Remark~\ref{rem:alt-fuchs}. Following
the ideas of A.~Khovanskii \cite{asik:fewnomials}, we introduce the
notion of a \emph{carpeting function} for \emph{noncompact}
semialgebraic subsets. Loosely speaking, the carpeting function
plays the role of the ``reciprocal distance to the frontier of the
set''.

\begin{Def} \label{def:carp-function}
A \emph{carpeting function} of a noncompact semialgebraic
subset $Z$ is a continuous semialgebraic positive function
$R\:Z\to\R_+$ which is proper as a map from $Z$ to $\R_+$. For
technical reasons we will always require in addition that $R\ge
2$ on $Z$.
\end{Def}

In other words, a continuous semialgebraic function $R$ is carpeting
for a non-compact set $Z$ if $R(z)$ tends to infinity along any
sequence $\{z_k\}_{k=1}^\infty$ without accumulation points in $Z$
(``converging to a frontier of $Z$'').

\begin{Ex}
The function $R_\flat$ is ``almost carpeting'' on $\clF$: it is
proper, yet not semialgebraic, since the Fubini--Study distance is
not algebraic. Yet the latter circumstance is purely technical: one
can easily construct a semialgebraic distance on $\C P^1$ and use it
in \eqref{eq:main-R} instead of the Fubini--Study distance.
\end{Ex}

Any two carpeting functions on the same set are related by a
two-sided \Loja.-type inequality.

\begin{Lem}\label{lem:lojasiewicz}
For any two carpeting functions $R_1,R_2$ on the same
semialgebraic set $Z$, there exist a finite positive
constant $s$ such that
\begin{equation}\label{eq:carp-equivalent}
    \forall z\in Z\qquad R_2^{1/s}(z)\le R_1(z)\le R_2^s(z).
\end{equation}
\end{Lem}

\begin{proof}
Consider the joint graph $\{(z,u,v)\in\R^{n+2}\:z\in
Z,\allowbreak u=R_1(z),\allowbreak v=R_2(z)\}$, which is a
semialgebraic set by construction. Its projection $S\subseteq
\R^2$ on the $(u,v)$-plane parallel to the $z$-direction is
semialgebraic by the Tarski--Seidenberg principle.

The function $\f(u)=\sup\{v\in\R_+\:(u,v)\in S\}$ takes
only finite values. Indeed, the continuous function $R_2$
attains its maximum on the compact subset $\{z\in
Z\:\R_1(z)=u\}$. Since the function $R_1$ is carpeting, it
assumes arbitrarily large values, hence the function $\f(u)$ is
defined for all sufficiently large $u\in[2,+\infty)$. Being
semialgebraic by construction, it grows no faster than
polynomially, $\f(u)\le C u^s$, $C,s\in\R_+$. The constant $C$
can be absorbed into the increased power, $\f(u)\le u^{s+c}$,
$c=\operatorname{log}_2 C$, for all $u\ge 2$.

The inequality in the other direction is obtained analogously.
\end{proof}

In the future any two functions $R_1,R_2$ constrained by the
inequalities \eqref{eq:carp-equivalent} on their common domain, will
be referred to as \emph{polynomially equivalent}.

\begin{Rem}\label{rem:alt-fuchs}
The class of Fuchsian systems (or, more precisely, its interior) can
be alternatively described as follows. Consider the rational matrix
function
\begin{equation}\label{exp-repr}
    A(t,\l)=\mathbf P(t)/Q(t),\qquad \mathbf P(t)=\sum _0^{m-1}
    M_k\,t^k,\quad Q(t)=t^m+\sum_0^{m-1}c_k\,t^k,
\end{equation}
with the $(n\times n)$-matrix coefficients $M_0,\dots,M_{m-1}$ and
scalar coefficients $c_0,\dots,\allowbreak c_{m-1}\in\C$ denoted by
$\l\in\C^d$, $d=(m-1)(n^2+1)$. Assume further that $\mathbf P(t)$
and $Q(t)$ are coprime.  The corresponding linear system will be
Fuchsian provided that the denominator has no multiple roots, that
is, the discriminant $\Delta\in\C[\l]$ of the polynomial $Q$ (an
explicit polynomial in $\l$) does not vanish (the point at infinity
is always Fuchsian).

There is an obvious carpeting function on the space of the
parameters,
\begin{equation}\label{R-sharp}
    R_\sharp(\l)=2+\frac1{|\Delta(\l)|}+\|\mathbf P\|+\|Q\|.
\end{equation}
One can easily establish a two-sided inequality between the
functions $R_\flat$ and $R_\sharp$  and show that they are
polynomially equivalent (one has to express the discriminant via
differences of the roots of the polynomial $Q$). Moreover, one can
easily show that the in upper bound \eqref{mainres} the function
$R_\flat$ can be replaced by $R_\sharp$ without changing the growth
rate of the exponent $\nu$, cf.~with Corollary~\ref{cor:any-carpet}
below.
\end{Rem}

\subsection{Product structure of the space $\clF$ of Fuchsian
 systems}\label{sec:prod-struct}
Representation \eqref{fuchs-principal-parts} together with the
special role played by the spectra of the residues and geometry of
the singular locus, suggests that the space $\clF$ of Fuchsian
systems is naturally (bijectively) parameterized by a semialgebraic
variety having a rather special \emph{product structure}:
\begin{equation}\label{eq:config-space}
\begin{aligned}
    \clF &\simeq\mathscr P_m\times \mathscr R_{n,m}
    \subseteq\bigl(\C P^1\bigr)^m\times\bigl(\C^{n^2}\bigr)^m,
    \\
    \mathscr P_m&=\{(\tau_1,\dots,\tau_m)\:\tau_i\in\C P^1,\
    \tau_i\ne \tau_j\text{ for }i\ne
    j\}\subseteq\bigl(\C P^1\bigr)^m,
    \\
    \mathscr R_{n,m}&=\{(A_1,\dots,A_m)\:A_i\ne 0,\ A_i\in\Mat_n(\C)\}
    \subseteq\bigl(\C^{n^2}\bigr)^m.
\end{aligned}
\end{equation}
This parametrization is continuous (even biholomorphic) in the sense
that convergence of the tuples $(\tau_i,A_i)\in\mathscr
P_m\times\mathscr R_{n,m}$ in the natural (geometric) sense implies
the uniform convergence of the corresponding rational matrix 1-forms
on $\C P^1$ on compact sets disjoint with the singular loci of the
forms (but not necessarily vice versa).

From now on we will identify the variety $\clF $ of Fuchsian matrix
1-forms on $\C P^1$ having $m$ distinct poles, with points
($2m$-tuples) from the space $\bigl(\C
P^1\bigr)^m\times\bigl(\C^{n^2}\ssm\{0\}\bigr)^m$. The set $\clS $
becomes then a relatively closed subset of $\clF$.

The compactification of the space $\clF $ preserving the
aforementioned  product structure is the product of projective
spaces:
\begin{equation}\label{eq:conf-space-comp}
    \overline{\mathscr F}_{n,m}=(\C P^1\bigr)^m\times\bigl(\C
    P^{n^2}\bigr)^m.
\end{equation}
The frontier $\partial\clF =\overline{\mathscr F}_{n,m}\ssm\clF$ of
the variety of Fuchsian systems consists of components of three
types,
\begin{equation}\label{eq:3-bds}
    \partial \clF \subseteq\mathscr V_{n,m}
    \cup\mathscr E_{n,m}\cup\mathscr C_{n,m},
\end{equation}
where:
 \begin{enumerate}
  \item the \emph{vanishing frontier} $\mathscr
V_{n,m}$ consists of the $2m$-tuples
 with one or more zero residue matrix
$A_i$ equal to zero,
 \item the \emph{explosion frontier} $\mathscr E_{n,m}$ is
the subset corresponding to one or more ``\emph{infinite}'' residues
$A_i$, and
 \item $\mathscr C_m=$ is the \emph{collision frontier}
corresponding to the union of the diagonals $\bigcup_{i\ne
j}\{\tau_i=\tau_j\}$.
 \end{enumerate}
The parametrization of Fuchsian 1-forms by poles and residues
extends continuously on the vanishing and the collision frontier,
but not on the explosion frontier. Note that while the passage to
limit as $\Omega$ tends to $\mathscr E_{n,m}$ is impossible, the
\emph{simultaneous limit} on the intersection $\mathscr
E_{n,m}\cap\mathscr C_{n,m}$ may well make sense in the class of
\emph{rational} systems. One can argue that generically such
``collision with exploding residues'' results in creation of an
\emph{irregular linear system}, thus the intersection $\mathscr
E_{n,m}\cap\mathscr C_{n,m}$ may well be considered as an
\emph{irregularity frontier.}

Examples~\ref{ex:norm-dependence} and~\ref{ex:distance-dependence}
suggest that the counting function $\cN\:\clS \to\R_+$ necessarily
grows to infinity near the explosion and the collision frontier
components. On the contrary, Theorem~\ref{thm:main} shows that the
counting function is bounded near the vanishing frontier.

Together with Theorem~\ref{thm:quasivk} this allows to argue that,
as long as the problem of counting zeros of solutions is considered
as an \emph{algorithmic mass problem formulated for the entire class
of general Fuchsian systems}, the natural domain for the counting
function should be the partial closure of the set $\clS $, obtained
by adjoining the vanishing frontier component $\mathscr V_{n,m}$,
\begin{equation}\label{eq:natdom}
\begin{aligned}
    \clS^*&=\{(\tau_1,\dots,\tau_m)\:\tau_i\ne \tau_j\}
    \times\{(A_1,\dots,A_m)\:\operatorname{Spec}A_i\subseteq\R\}
    \\
    &\subseteq\clF^*=\mathscr P_m\times\bigl(\C^{n^2}\bigr)^m.
\end{aligned}
\end{equation}
The partial closures $\clF^*$ and $\clS^*$ are semialgebraic
varieties ``parameterizing'' the class of Fuchsian systems with $m$
singularities and eventually zero residue matrices (resp., the class
of such systems satisfying the spectral condition \eqref{eq:spec}).

The function $R_\flat$ defined by the formula \eqref{eq:main-R}, is
clearly a positive function on $\clS ^*$ which tends to infinity
polynomially as $\Omega$ approaches the frontier of $\clF $. This
function is not semialgebraic because the Fubini--Study distance
$\dist\:\C P^1\times\C P^1\to\R_+$ is not semialgebraic, yet this
failure is purely technical; clearly, there are semialgebraic
distance functions on the projective line, all of them equivalent to
each other; this would transform \eqref{eq:main-R} into a genuine
carpeting function. However, this choice is by no means unique.

Since all carpeting functions on the same semialgebraic set are
equivalent, Theorem~\ref{thm:main} admits the following
reformulation.

\begin{Cor}\label{cor:any-carpet}
For any carpeting function $R\:\clS ^*\to\R_+$ the counting function
$\cN(\Omega)$ admits an explicit polynomial upper bound of the form
$$
 \cN(\Omega)|_{\clS ^*}\le R^\nu(\Omega)\qquad\text{for
 some finite }\nu.
$$
The constant $\nu$ depends on $n,m$ and the choice of $R$.
\end{Cor}

\subsection{Product spaces and natural carpeting functions on them}
There is no single distinguished (or preferred) carpeting
function on a semialgebraic variety, and it is even less clear how one may
compare carpeting functions defined on different semialgebraic
varieties. Thus it is rather difficult to analyze the asymptotic
behavior (in $n,m$) of different bounds in invariant terms.

However, in the particular case of the varieties of Fuchsian systems
$\clF $ one may argue that symmetry considerations distinguish a
certain class of \emph{natural} carpeting functions. More precisely,
we introduce two special classes of \emph{product semialgebraic
spaces}, provisionally called \emph{Fermi and Bose types}
(respectively, with or without the explicit prohibition of
coincident terms in the product) of several copies of the same
semialgebraic variety $Z$. Then we show that on such product spaces
one can introduce ``anonymous'' carpeting functions that do not
depend \emph{explicitly} (in the accurate sense introduced below) on the number of
the copies. These natural carpeting functions depend on the product
structure types and are defined modulo \emph{parametric polynomial
equivalence} (as described in
Definition~\ref{def:carp-families-equivalent}).

Let $Z$ be an arbitrary semialgebraic set.

\begin{Def}
A \emph{Bose product space} $Z^n$ is the (Cartesian) power of
$Z$, the space of all tuples of points $\(z_1,\dots,z_n\)$,
$z_i\in Z$.

A \emph{Fermi product} $Z^{*n}$ is the space of all \emph{pairwise
different} tuples of points $\(z_1,\dots,z_n\)$, $z_i\in Z$:
\begin{equation*}
    Z^{*n}=Z^n\ssm\bigcup_{i\ne j}\{z_i=z_j\}.
\end{equation*}
\end{Def}

Assume now that $R\:Z\to\R_+$ is a carpeting functions on $Z$ and
$R'$ is a symmetric carpeting function on the Fermi square $Z^{*2}$,
$R'(z_1,z_2)=R'(z_2,z_1)$. Starting from these two functions, we can
construct carpeting functions on arbitrary Bose/Fermi products,
aggregating these ``basic'' functions by an arbitrary \emph{semialgebraic
operation}.

\begin{Def}\label{def:carp-op}
A continuous commutative and associative binary operation
$\odot\:\R_+\times\R_+$ will be called \emph{semialgebraic}, if its
graph is a semialgebraic subset in $\R^3_+$, and \emph{carpeting},
if all sets $\{(x,y)\in\R^2_+\:x\odot y\le C\}$ are compact and
exhaust $\R_+\times\R_+$ as $C\to+\infty$ (i.e., if the function
$f(x,y)=(x\odot y)+2$ is carpeting in the sense of
Definition~\ref{def:carp-function}).
\end{Def}

Because of the commutativity and associativity, the expressions
$\bigodot_{i=1}^n c_i$ make sense for any unordered collection of
positive numbers $c_i$.

\begin{Ex}
The sum $x+y$, maximum $\max(x,y)$, ``radius'' $\sqrt{x^2+y^2}$
are semialgebraic carpeting operations. The product is
semialgebraic albeit not carpeting, yet the ``shifted product''
$x\odot y=(2+x)(2+y)$ is.

More generally, let $h\:\R_+\to\R_+$ be any semialgebraic
homeomorphism of $\R_+$ and $\odot$ any semialgebraic carpeting
binary operation (e.g., one from the above list). Then the
binary operation $(x,y)\mapsto h^{-1}\bigl(h(x)\odot
h(y)\bigr)$ is again carpeting.
\end{Ex}

\begin{Def}\label{def:natural-carp}
A \emph{natural} family of carpeting functions on the Bose products
$Z^n$ (resp., Fermi products $Z^{*n}$) is a family of carpeting
functions defined by the formulas
\begin{equation}\label{bigodot}
    R(z_1,\dots,z_n)=\bigodot_{i=1}^n
    R(z_i),\quad \text{resp.,}\quad
    R'(z_1,\dots,z_n)=\bigodot_{i\ne j}R'(z_i,z_j),
\end{equation}
where $R$ (resp., $R'$) is an arbitrary semialgebraic carpeting
function on $Z^1$ (resp., on $Z^{*2}$) and $\odot$ any semialgebraic
carpeting commutative associative binary operation on the
nonnegative ray $\R_+$.
\end{Def}

Combining arbitrary \emph{natural} carpeting functions for the Bose
powers $\mathscr R_{n,m}\simeq \bigl(\C^{n^2}\bigr)^m$ and Fermi
powers $\mathscr P_m=(\C P^1)^{*m}$, by the carpeting operation
$\odot$, we obtain the class of \emph{natural carpeting functions}
on the product spaces $\clF =\mathscr P_m\times\mathscr R_{n,m}$ as
in \eqref{eq:config-space}.

\subsection{Polynomial equivalence of natural carpeting
 functions}
The degree of freedom used to introduce the natural carpeting
functions on the product spaces is considerable, as one can vary
both the basic functions $R,R'$ as well as the semialgebraic
operation $\odot$. However, these variations result in controllable
change of the growth rate.

\begin{Def} \label{def:carp-families-equivalent}
Two \emph{families} of carpeting functions $R_1,R_2\:Z_n\to\R_+$
defined on the same family of semialgebraic sets $Z_n$ indexed by a
natural parameter $n$ (in general, the dimensions $\dim Z_n$ grow to infinity with $n$), are said to be \emph{polynomially equivalent}, if there
exists a sequence of finite positive constants $s_n\in\R_+$, such
that,
\begin{equation}\label{eq:carp-family-eq1}
    R_2^{1/s_n}(z)\le R_1(z)\le R_2^{s_n}(z)\qquad \forall z\in Z_n;
\end{equation}
cf.~with \eqref{eq:carp-equivalent}, and this sequence grows at most
polynomially as $n\to\infty$:
\begin{equation}\label{eq:carp-family-eq2}
    \exists c,r<+\infty \quad \text{ such that }\quad s_n\le cn^r,\qquad
    \forall n=1,2,3,\dots
\end{equation}
An obvious modification allows one to also speak of polynomially
equivalent carpeting functions for families $Z_\alpha$ indexed by a
multiindex $\alpha=(\alpha_1,\dots,\alpha_n)$. In this case the
respective exponents $s_\alpha\in\R_+$ should grow no faster than
$c\,|\alpha|^r$ for some finite $c,r\in\R$.
\end{Def}

Clearly, there is no reason why two arbitrary families of carpeting
functions be polynomially equivalent: even if this happens, one can
involve one of the functions into a sufficiently quickly growing
sequence of powers and destroy the initial equivalence. Yet any two
\emph{natural} parametric families of carpeting functions on a
family of product spaces turn out to be \emph{always equivalent} in
the sense of Definition~\ref{def:carp-families-equivalent}. This can
be considered as a parametric version of the \Loja. inequality.

In complete analogy with Lemma~\ref{lem:lojasiewicz}, one can show
that any two carpeting operations are equivalent, in particular, for
any such operation $\odot$,
\begin{equation}\label{eq:carp-loja}
   \exists\gamma<+\infty,\ \forall a,b\ge 2,
   \qquad  (ab)^{1/\gamma}\le a\odot b\le (ab)^\gamma.
\end{equation}
Iterating this inequality, we obtain the following result.

\begin{Lem}\label{lem:poly-equiv}
For any carpeting operation $\odot$ as in
Definition~\ref{def:carp-op}, there exist two finite positive
constants $c,r$ such that for any $n\in\N$ and any nonnegative
numbers $x_1,\dots,x_n\ge 2$
\begin{equation}\label{eq:carp-iter}
    x_1\odot\cdots\odot x_n\le (x_1\cdots x_n)^{cn^r}.
\end{equation}
\end{Lem}

\begin{proof}
It is sufficient to prove the inequality for the number of terms
being a power of two (for the intermediate cases one can replace
some of the extra terms by units).

For $n=2^m$ the proof goes by induction, iterating
\eqref{eq:carp-loja}. Indeed, if $\gamma<2^r$, then  the
$\odot$-product of $2^{m+1}$ terms is bounded by the product as
follows,
 $$
 (X_1^{c\,2^{mr}}X_2^{c\,2^{mr}})^\gamma\le (X_1X_2)^{\gamma
 c\,2^{mr}}\le (X_1X_2)^{c\,2^{(m+1)r}},
 $$
where $X_1=x_1\cdots x_{2^m}$ and $X_2=x_{2^m+1}\cdots
x_{2^{(m+1)}}$ stand for the products of the first and the last
$2^m$ terms. The base of the induction can always be guaranteed by a
sufficiently large $c>0$.
\end{proof}

\begin{Lem}\label{lem:nat-poly-eq}
Any two natural carpeting families on Bose \textup(resp.,
Fermi\textup) powers $Z^n$ \textup(resp., $Z^*n$\textup) of the
same semialgebraic set $Z=Z^1$ are polynomially equivalent to
each other.
\end{Lem}

\begin{proof}
The proof follows immediately from Lemma~\ref{lem:lojasiewicz} and
Lemma~\ref{lem:poly-equiv}.
\end{proof}

The function $R_\flat$ defined by formula \eqref{eq:main-R},
represents the class of natural carpeting functions on the product
spaces $\clF^*$ and its relatively closed subset $\clS^*$. From
Lemma~\ref{lem:nat-poly-eq} we conclude that the assertion of
Theorem~\ref{thm:main} does not depend on the specific choice of the
expression \eqref{eq:main-R}; the counting function $\cN(\cdot)$ for
zeros of solutions of Fuchsian systems admits semialgebraic bound
which is double exponential in the complex dimension $d=\dim
\clF^*\sim n^2m$ of the corresponding variety of Fuchsian systems.

The following result is the most general form of the polynomial
bound on zeros of solutions of Fuchsian systems.

\begin{Thm}\label{thm:main-universal}
For any natural carpeting function $R(\cdot)$ on the class of
Fuchsian systems considered as a product space $\clF^*$,
\begin{equation}\label{ex:2exp-univ}
 \mathcal N\,\big|_{\clS^*}\le R^{2^{O(d)}},
 \qquad d=\dim_{\C}\clF^*.
\end{equation}
The constant in $O(d)$ is explicit, though it depends on the
choice of the natural carpeting function.
\end{Thm}

\begin{proof}
This follows immediately from Theorem~\ref{thm:main} since the
carpeting functions \eqref{eq:main-R} form a natural carpeting
family, and the polynomially growing terms $s_{n,m}$ of the
polynomial equivalence are obviously absorbed by the
exponential $2^{O(m^2n)}$.
\end{proof}

\begin{Rem}
Note that the carpeting function $R_\sharp$ introduced in
\eqref{R-sharp}, is \emph{not natural} on the space $\clF$ with
respect to the product structure of the latter as it was introduced
in Definition~\ref{def:natural-carp}. Despite that, the function
$R_\sharp$ is polynomially equivalent to the function $R_\flat$ in
the sense of Definition~\ref{def:carp-families-equivalent}.
\end{Rem}

\subsection{Near optimality of the double exponential bounds}
The proof of Theorem~\ref{thm:main-universal} suggests that the
double exponential bound asserted in \eqref{ex:2exp-univ} is
\emph{almost optimal}. Indeed, a more moderate growth rate (say,
with the exponent $\nu$ in \eqref{mainres} growing polynomially in
$n,m$, i.e., just one notch slower than the exponent) would already
be highly dependent on the choice of carpeting function. This
observation suggests that, at least \emph{as long as the entire
Fuchsian class $\clS^*$ is concerned}, one cannot expect very
significant improvements to the (non-uniform) upper bounds of
Theorem~\ref{thm:main-universal}.

Of course, the situation changes dramatically if instead of the
entire class $\clF$ we consider its various subclasses. In
\cite{inf16} we show that for a \emph{isomonodromic regular
quasiunipotent classes} of linear systems on $\C P^1$ one can
produce a uniform upper bound for the number of zeros of isolated
solutions. This bound also turns out to be double  exponential, yet
because of the uniformity one can easily ask about its possible
improvement.

The rest of the paper is devoted to the detailed proof of
Theorem~\ref{thm:main}.

\section{From Fuchsian systems to Fuchsian equations with a complexity control}

It is now well understood that there is a substantial difference between
systems of first order linear ordinary differential equations and
(scalar) high order linear ODE's with respect to studying zeros of
solutions and their linear combinations, see
\cite{mit:counterexample} and references therein. In this section we
prove that any linear combination of solutions of a Fuchsian system
\eqref{eq:ls} from the class $\clF^*$ satisfies a linear equation of
order $n$ with \emph{explicitly bounded coefficients}.

\begin{Lem}\label{lem:scalar-eq}
An arbitrary linear combination of coordinates of a solution of a
Fuchsian system with the matrix 1-form $\Omega\in\clF^*$, in a
suitable affine chart $t$ on $\C P^1$, satisfies a linear ordinary
differential equation with polynomial coefficients,
\begin{equation}\label{eq:le}
\begin{gathered}
    a_0(t)\,y^{(n)}+\cdots+a_{n-1}(t)\,y'+a_n(t)\,y=0,
    \\
    a_0,\dots,a_n\in\C[t],\qquad\gcd(a_0,\cdots,a_n)=1
\end{gathered}
\end{equation}
and the following qualifications\textup:
\begin{enumerate}
 \item [(a)]$\deg a_j\le d=n^2m$,
 \item [(b)]$\norm{a_j}\le \norm{a_0}\cdot R^\nu$, with
 $R=R_\flat(\Omega)$ and $\nu\le 2^{O(d)}$,
 \item[(c)] solutions of the equation \eqref{eq:le} are holomorphic outside the
 finite locus $\S=\{\tau_1,\dots,\tau_m\}\subset\C$ such that
 \begin{equation}\label{eq:locus}
  |\tau_i|\le R,\qquad |\tau_i-\tau_j|\ge 1,
 \end{equation}
 in particular, solutions are holomorphic in the disk $\{|t|>R\}$.
\end{enumerate}
\end{Lem}

The proof of this Lemma is largely parallel to the exposition from
\cite[\parasymbol2]{ongrig}. For convenience of the reader we
reproduce it with the required modifications. In what follows we will
refer to the ratio
\begin{equation}\label{pre-slope}
    \slope D=\max_{j=1,\dots,n}\frac{\|a_j\|}{\|a_0\|}
\end{equation}
as the \emph{slope} of the differential operator $D$ as in the left
hand side \eqref{eq:le} (and the corresponding linear ordinary
homogeneous equation). Thus Lemma~\ref{lem:scalar-eq} asserts an
explicit global upper bound on the slope of the equation
\eqref{eq:le} derived from the Pfaffian system $\d x=\Omega x$, in
terms of the carpeting function $R_\flat(\Omega)$.

\subsection{Preparatory remarks}
We start with choosing a convenient affine chart. From elementary
geometric arguments involving spherical areas, it follows that for
any $m$ points on the Riemann sphere, there always exists a point at
least $O(1/\sqrt m)$-distant from them in the sense of Fubini--Study
distance. We choose the affine chart $z\in\C$ in which this point
corresponds to infinity and normalize it so that the Fubini--Study
metric is given by the form $\frac{|dz|}{1+|z|^2}$. Such map is
defined uniquely modulo the Euclidean rotation, and we have an
obvious inequality $\dist(z_1,z_2)\le |z_1-z_2|$.

By construction, in the chosen chart $z$ the affine coordinates
$\tau_1,\dots,\tau_m$ of the singular points of a system \eqref{eq:ls} with
$R_\flat(\Omega)=R$, satisfy the inequality $|z_j|\le O(\sqrt m)$
for all $j=1,\dots,m$. Since all terms in the expression
\eqref{eq:main-R} are nonnegative and the residues are invariant, we
conclude that in this chart the inequalities
\begin{equation}\label{eq:prep-bounds}
    |z_i|\le O(\sqrt m),\quad j=1,\dots,m,
    \qquad |z_i-z_j|\ge 1/R,\quad i\ne j
\end{equation}
hold simultaneously. By a suitable affine rescaling $t=cz$ we can
normalize the singular locus of the system to satisfy the
inequalities \eqref{eq:locus}. In the affine chart thus constructed,
the Pfaffian linear system \eqref{eq:ls} is equivalent to a system
of linear ordinary differential equations of the form
\begin{equation}\label{eq:lst}
\begin{gathered}
    \frac{\d x}{\d t}=A(t)\cdot x,\quad
    A(t)=\sum_{j=1}^m\frac{A_j}{t-\tau_j},\\
    |A_j|\le R,\quad |\tau_j|\le R,\qquad |\tau_i-\tau_j|\ge 1.
\end{gathered}
\end{equation}

Finally, performing a linear transformation $x\mapsto Ux$ with a
unitary (constant) matrix $U$, we may without loss of generality
assume that the linear combination in question is just the first
component, $y(t)=x_1(t)$. This transformation does not affect the
norms of the residues, hence the inequalities \eqref{eq:prep-bounds}
are preserved.

\subsection{Objects defined over $\Q$ and their
 complexity}\label{sec:complexity}
If we consider the entries of the residue matrices $A_i$ and the
affine coordinates $\tau_i\in\C$ of the singular points as
coordinates on $\clF^*$, denoting a general point of this space by
$\l$, then the coefficients of the linear system \eqref{eq:lst}
become \emph{integer numbers}, more precisely, only $0$ and $\pm 1$.
In other words, the entries of the matrix 1-form $\Omega_\l$ are
defined over the field $\Q$ on the space with coordinates $(t,\l)$.
In sharp contrast with the fields of real or complex numbers, for
objects defined over the field $\Q$ one can define their
complexity (bitsize, or, more precisely, the exponent of the
bitsize). This definition may be introduced by many equivalent ways. For instance, as follows:\begin{enumerate}
 \item Complexity of a polynomial $p\in\Z[\l,t]$ is its norm
 $\|p\|\in\Z_+$, the sum of all absolute values of its integer
 coefficients;
 \item Complexity of a rational function $R\in\Q(\l,t)$ is the minimal value of
 sum  $\|p\|+\|q\|$ over all representations $R=p/q$ with
 $p,q\in\Z[\l,t]$; in particular, complexity of an irreducible
 rational number $r=p/q$, $p,q\in\Z$, $\gcd(p,q)=1$, is $|p|+|q|$;
 \item Complexity of any object constructible over $\Q(\l,t)$ is the
 sum of complexities of all rational functions explicitly occurring
 in the construction. For instance, complexity of a rational matrix
 1-form $\Omega$ is the sum of complexities of all its entries;
 complexity of a differential operator with rational coefficients
 is the sum of complexities of its coefficients \etc.
\end{enumerate}
This definition (as most constructions from the complexity theory)
in fact applies not to the object in question, but rather to a
specific \emph{representation} of this object. It will be important
in this paper that the  objects under discussion are described by
their \emph{explicit} representation.

Whenever two representations are related by a simple transformation
(e.g., expansion or reduction to common denominator), their
complexity is usually comparable. More generally, \emph{any explicit
construction over $\Q$ admits explicit control over the complexity
growth of all intermediate objects}. This observation is the primary
reason for introducing the objects over $\Q$ in the context where
the initial field is $\C$ or $\R$. We will apply this technique to
the straightforward process of derivation of a scalar differential
equation of high order from a system of first order equations
(consecutive derivations and elimination).

\subsection{Derivation of the scalar Fuchsian equation over $\mathbb Q(\l,t)$}
 \label{sec:eq-derivation}
By direct inspection of the formula \eqref{eq:lst} one can see that
the matrix product $\mathbf P=\prod_1^m (t-\tau_i)A(t)$ is a polynomial
matrix with entries in $\Z[\l,t]$ of degree $\le m$ and norm not
exceeding $m\cdot2^{m-1}\le 2^m$ (here and below we use the fact
that the norm is multiplicative, $\|pq\|\le \|p\|\cdot\|q\|$ for any
two polynomials $p,q$). For brevity, we call elements
of the ring $\Z[\l,t]$ the \emph{lattice polynomials}.

\begin{Lem}\label{lem:deriv}
The first component of any solution of the system \eqref{eq:ls}
considered as a linear system with coefficients from the field
$\mathbb Q(\l,t)$, satisfies a linear $n$th order differential
equation with lattice polynomial coefficients of the form
\begin{equation} \label{eq:linear}
\begin{gathered}
   a_{0}(\l,t)\,y^{(n)}+a_{1}(\l,t)\,y^{(n-1)}+\cdots+a_{n}(\l,t)\,y=0,
   \\
   a_0,\dots,a_n\in\Z[\l,t],\qquad a_0\ne0.
\end{gathered}
\end{equation}

The coefficients $a_j$ are lattice polynomials in the variables
$\l,t$ having degrees $\le n^{2}m$ and explicitly bounded norms,
$\norm{a_j}\le 2^{O(n^2m)}$.
\end{Lem}

Note that the norms in the formulation of Lemma~\ref{lem:deriv} are
computed with respect to the ring $\Z[\l,t]$, hence they are
\emph{nonnegative integer numbers}.

\begin{proof}[Proof of the Lemma]
By virtue of the system \eqref{eq:lst}, the higher derivatives of
the function $y(t)$ can be expressed inductively as linear
combinations of the dependent variables $x=(x_1,\dots,x_n)$,
\begin{equation} \label{eq:higher-derivatives}
   y^{(k)}(t)=\xi_{k}\cdot x,\qquad \xi_k\in\bigl(\mathbb
   Q(\l,t)\bigr)^m,
\end{equation}
with the $n$-vectors $\xi_k$ over the field $\mathbb Q(\l,t)$
determined by the recursive rule
\begin{equation} \label{eq:coefficients-definition}
\begin{aligned}
   \xi_{0} &= \(1,0,\ldots,0\), \\
   \xi_{k+1}&= \partial_t \xi_{k}+
     \xi_{k}\cdot A,
\end{aligned}
\qquad A\in\Mat_n\bigl(\mathbb Q(\l,t)\bigr), \
\partial_t=\partial/\partial t
\end{equation}
(we treat $x$ as a column vector and $\xi_k$ as row vectors from
the dual space).

The structure and complexity of the $n$-covectors $\xi_k$ can easily
be established by the inspection of the formulas
\eqref{eq:coefficients-definition}. Denote
\begin{equation*}
   Q(\l,t)=\prod_{i=1}^{m}(t-\tau_{i})\in\Z[\l,t], \qquad
   \|Q\|=2^m,\ \deg Q=m,
\end{equation*}
the \emph{discriminant}, so that $A=\mathbf P/Q$ in
\eqref{eq:coefficients-definition}.

The $n$-vectors $\xi_{0},\ldots,\xi_{n}$ defined by
\eqref{eq:coefficients-definition}, satisfy the following
properties:
\begin{enumerate}
\item Their multiples $\alpha_k=Q^k\xi_k$ are lattice $n$-vector
polynomials for all $k=0,\dots,n$, of degree no greater than $km$,
and the norms $\norm{\alpha _{k}}\le (3nm)^{n}2^{nm}$.
\item The wedge product of any $n$ vectors from this collection, the
polynomial $n$-form
\begin{equation*}
   Q^{n(n+1)/2}\xi_{0}\wedge\cdots \wedge {\xi _{i-1}}\wedge {\xi _{i+1}}\wedge \cdots \wedge
   \xi_{n},
\end{equation*}
has as its only coefficient a lattice polynomial of degree at most
$n^{2}m$ and norm bounded by $n!(3nm)^{n^{2}}2^{2n^{2}m}$.

\item The wedge product $\xi_{0}\wedge \cdots \wedge\xi_{n-1}$ does not
vanish identically on $\clF^*$.
\end{enumerate}

Indeed, the polynomial 1-forms $\alpha_{k}=Q^k\xi_k$ satisfy the
recurrence relations
\begin{equation}
      \alpha_{k+1}=Q\cdot \partial_t \alpha_{k}-(k\,\partial_t Q-\mathbf P)\cdot
      \alpha_{k},
\end{equation}
which follow directly from \eqref{eq:coefficients-definition}. The
bound for the degrees $\deg\alpha_k$ is obvious since $\deg Q,\deg
\mathbf P\le m$. The bound for the norms $\|\alpha_k\|$ can also be
derived from the recursive equations,
\begin{eqnarray*}
   \norm{\alpha_{k+1}} &=& \norm{Q\cdot \partial_t \alpha_{k}-(k\partial_t Q-\mathbf P)\alpha_{k}}\\
   &\le & \norm Q (\deg \alpha _{k}) \norm{\alpha_{k}}+k(\deg Q)\norm Q
      \norm{\alpha_{k}} +n\norm P \norm {\alpha_{k}} \\
   &\le & (2^{m}km+2^{m}km+2^{m}nm) \norm{\alpha_{k}} \\
   &\le & 3\cdot 2^{m}\cdot nm \norm{\alpha_k},
\end{eqnarray*}
so that an easy induction gives $\norm{\alpha _{k}}\le
(3nm)^{n}2^{nm}$. To establish the bound for the wedge products, we
expand the $n\times n$ determinants using the Laplace expansion
involving $n!$ summands, each having the degree $\le
\frac12{n(n-1)m}$ and norm $\le (3nm)^{n^{2}}2^{n^{2}m}\le2^{O(n^2m)}$.

Finally, to establish (3), we note that for \emph{specific}
Fuchsian systems \eqref{eq:ls} and for specific solutions
$y(t)$ we can apriori guarantee that the functions
$y^{(0)}(t),\ldots,y^{(n-1)}(t)$ are linearly independent. For
instance, the function $y(t)=1+t+\cdots+t^{n-1}$ is obviously
linearly independent with its derivatives of all orders $\le
n-1$. This function is a linear combination of solutions of an
Euler system with two singularities at the origin and at
infinity, which is a particular case of a Fuchsian system. For
the corresponding value of the parameters $\l$, the Wronskian
$\xi_0\wedge\cdots\wedge\xi_{n-1}$ is not identically zero as a
function of $t$.

To derive the equation \eqref{eq:linear}, it remains to note that
any $n+1$ vectors in the $n$-space $\bigl(\mathbb Q(\l,t)\bigr)^n$
are linearly dependent over this field. In our particular case the
vectors $\xi_0,\dots,\xi_{n-1}$ are linearly independent, so the
vector $\xi_n$ expands as a linear combination,
$\xi_n=c_0\xi_0+\cdots+c_{n-1}\xi_{n-1}$, of the previous vectors
with coefficients $c_i\in \mathbb Q(\l,t)$. These coefficients can
be found using the Cramer rule, as ratios of the determinants,
\begin{equation}\label{eq:ratios-dets}
    c_i=\frac{\xi_{0}\wedge\cdots \wedge
    {\xi _{i-1}}\wedge\xi_{i+1}\wedge \cdots \wedge
   \xi_{n}}{\xi_{0}\wedge\cdots \wedge  \xi_{n-1}},
   \qquad i=0,\dots,n-1.
\end{equation}
The numerator and the denominator of any such fraction are
rational functions defined over $\Q$ with the denominators being known powers of the
discriminant $Q$. Multiplying the linear identity by the common
denominator, we obtain a linear identity between the $n$-vectors
$\xi_0,\dots,\xi_n$ with \emph{lattice} polynomial coefficients, which by
\eqref{eq:higher-derivatives} produces a linear $n$th order
differential equation for the function $y$.

The leading coefficient $a_0$ of this equation is a nonzero lattice
polynomial, therefore its norm is a positive integer number which is
necessarily $\ge 1$.
\end{proof}

\subsection{Explicit bounds on coefficients
 of the scalar equation over $\C(t)$}
  \label{sec:eq-behaviour}
Specializing the equation \eqref{eq:linear} for each particular
value of the parameters $\l$ would produce an equation of the type
\eqref{eq:le} with coefficients from $\C[t]$. The bound for the
degrees of these coefficients will obviously be the same as for
\eqref{eq:linear}. However, the condition $\norm{a_0}_{\Z[\l,t]}=1$
does not prevent the specialization $a_0(\l,\cdot)$ of the leading
coefficient from having an arbitrarily small norm
$\norm{a_0(\l,\cdot)}_{\C[t]}$ or simply vanishing. In this
subsection we prove, following \cite{alexg:arxiv,ongrig}, that this
vanishing may happen if and only if all other (nonleading)
coefficients also vanish. Moreover, we show that the ratios
$\norm{a_j(\l,\cdot)}_{\C[t]}/\norm{a_0(\l,\cdot)}_{\C[t]}$ are
explicitly bounded as functions of $\l$. (Lemma~\ref{lem:explicit} below).

Consider a linear parametric differential equation \eqref{eq:linear}
with lattice polynomial coefficients $a_0,\dots,a_n\in\Z[t,\l]$.
Denote by $\varLambda\subset\C^\nu$ the singular locus corresponding
to the parameters $\l$ for which the leading coefficient
$a_0(\l,\cdot)$ vanishes identically as a univariate polynomial.

Solutions of the equation \eqref{eq:linear} may be considered as
holomorphic functions of $t,\l$ by fixing an initial condition at
some point $t=t_0$ and treating $\l$ as parameters of the equation
varying near some value $\l_0$, provided that the bidisk
$\{|t-t_0|<a,\ |\l-\l_0|<b\}$ is free from singular points of the
equation. Yet in general these solutions will exhibit singularities
on the locus $\{(\l,t)\:a_0(\l,t)=0\}\subset\C\times\C^\nu$ which
contains some exceptional lines $\{\l_*\}\times\C$ for
$\l_*\in\varLambda$, while intersecting the other (generic) lines
$\{\l\}\times\C$, $\l\notin\varLambda$, only by isolated points.

\begin{Def}\label{def:apparent-sing}
The family of equations \eqref{eq:linear} has an \emph{apparent
singularity} at a point $\l_*\in \varLambda\subset\C^n$, if there
exist:
\begin{enumerate}
 \item a small neighborhood $U\subset\C^\nu$ of $\l_*$,
 \item an open disk $D\subset\C$,
 \item a fundamental system of solutions $y_1(\l,t),\dots,y_n(\l,t)$ of the family,
\end{enumerate}
such that each function $y_i(\l,\cdot)$  analytically depends on
$\l$ in the disk $D$ as long as $\l\in U\ssm\varLambda$, and remains
bounded in $D$ as $\l\to\l_*$.
\end{Def}

By the removable singularity theorem, this means that a fundamental
system of solutions admits uniform analytic limit in $D$ as
$\l\to\l_*$, though the linear independence of the functions
$y_i(\l_*,\cdot)$ may fail at the limit.

Denote by $D_\l$ the linear differential operator corresponding to
the linear homogeneous equation \eqref{eq:linear} after specializing
the value $\l\notin\varLambda$. Denote by
\begin{equation}\label{slope}
    \f(\l)=\slope
    D_\l=\max_{j=1,\dots,n}\frac{\|a_j(\l,\cdot)\|}{\|a_0(\l,\cdot)\|},
    \qquad\l\in\C^\nu\ssm\varLambda,
\end{equation}
the slope of this operator. Without loss of generality we may assume
that the polynomials $a_0(\l,\cdot),\dots,a_n(\l,\cdot)$ are
mutually prime in the ring $\C[t]$: nontrivial common divisor may
occur only on a proper algebraic subset of $\C^\nu$ which can be
combined with $\varLambda$.

\begin{Lem}\label{lem:explicit}
Assume that the family of equations \eqref{eq:linear} with lattice
polynomial coefficients of known degree $\le d$ and bounded norms
$\norm{a_j}\le M$, depending on the complex $\nu$-dimensional
parameter $\l\in\C^\nu$, has only apparent singularities at all
points of the locus $\varLambda$ which is a proper algebraic subset of $\C^\nu$. Then
\begin{equation} \label{eq-ratio-bnd}
   \forall\l\text{ such that}\quad |\l|\le M,\qquad \f(\l)\le M^{d^{O(\nu)}}
\end{equation}
with an explicit constant in $O(\cdot)$.
\end{Lem}

\begin{proof}[Sketch of the proof]
The proof of the Lemma almost literally reproduces that of Lemma~5
from \cite{ongrig} with minimal modifications.

More specifically, we start with the observation that the slope $\f$
is locally bounded everywhere on $\C^\nu\ssm\varLambda$. This
follows from the assumption that the family \eqref{eq:linear} has
only apparent singularities, in the same way as in the proof of
Lemma~4 from \cite{ongrig}. The only required change is to replace
the unit disk by the disk $D$ from the
Definition~\ref{def:apparent-sing}.

Next, we consider the positive subgraph $X$ of the function
$\f(\cdot)$ over the ball of radius $M$. This subgraph is a
semialgebraic subset in $X\subseteq\C^{\nu+1}$, defined by lattice
polynomial equalities and inequalities. Since the function $\f$ is
locally bounded near each point of the compact ball $\{|\l|\le M\}$,
the semialgebraic set $X$ is globally bounded.

The lattice polynomials defining $X$ are of degrees not exceeding
$O(d)$ in $\Re\l$ and $\Im \l$. In the same way the norms of these
polynomials are no greater than $O(M)$. Finally, the dimension of
the space in which $X$ resides, is $O(\nu)$. By the quantitative
version of the Tarski--Seidenberg quantifier elimination theorem
\cite{grigoriev-vorobjov,heintz-roy-solerno} and, most recently,
\cite[Theorem~7.2]{basu-vorobjov}, the diameter of $X$ does not
exceed $M^{d^{O(\nu)}}$, as asserted in  \eqref{eq-ratio-bnd}.
\end{proof}

\begin{proof}[Proof of Lemma~\ref{lem:scalar-eq}]
The proof directly follows from Lemmas~\ref{lem:deriv}
and~\ref{lem:explicit}: we first derive the scalar equation
\eqref{eq:linear} over $\mathbb Q(\l,t)$ and transform it to the
form with lattice polynomial coefficients by eliminating the
denominators.

The resulting family of scalar equations \eqref{eq:linear} has only
apparent singularities, since for any parameter $\l$ solutions of
the initial system \eqref{eq:lst} depend analytically on $\l$ in any
disk free from singular points $\S=\S_\l\subseteq\C$ of the latter.
By Lemma~\ref{lem:explicit}, we conclude with the inequality
asserted by the Lemma, for all values $\l$ outside from the
vanishing locus $\varLambda=\{\l\:a_0(\l,\cdot)\equiv0\}$. This
locus is a proper algebraic hypersurface of the complex affine
space, hence by continuity the inequality
$\norm{a_j}\le\norm{a_0}\cdot R^\nu$ holds everywhere on the ball
$\{|\l|\le R\}$.
\end{proof}

\section{Counting zeros of solutions of Fuchsian equations} \label{sec:counting-zeros}

In this section we complete the proof of the main
Theorem~\ref{thm:main}. It is based on a complex counterpart of the
de la Vall\'ee Poussin theorem \cite{fields}, which asserts that the
variation of argument of any (complex) solution of a \emph{monic}
linear differential equation (with the leading coefficient
identically equal to $1$) along any rectilinear segment can be
explicitly bounded from above in terms of the magnitude of the
non-principal coefficients of this equation.

\subsection{Lower bounds of polynomials away from their zeros}
We start with an observation which relates the norm
\eqref{eq:poly-norm} on univariate polynomials with the
supremum-norm on disks of finite radius centered at the origin. Clearly, for any finite $d$
the two norms on the finite-dimensional space of polynomials of
degree $\le d$ are equivalent.  However, we shall also be interested
in the explicit constants involved in this equivalence.

The first side of this equivalence is easy to describe explicitly. Indeed,
$\|p\|=1$ implies an easy upper bound $|p(t)|\le R^d$ for $|t|\le R$. We
will prove the inequality in the opposite direction, establishing
a \emph{lower} bound for $|p(t)|$ for all $t\in\C$ distant from the null
locus of $p$, provided that $\|p\|_{\C[t]}=1$.

\begin{Lem}\label{lem:lower-bound}
Let $p=\sum_0^d p_j\,t^j$ be a polynomial of degree $d$ and unit
norm $\|p\|_{\C[t]}=1$, and $\S=\{t\in\C\:p(t)=0\}$ the zero locus
of $p$.

Then for an arbitrary point $t\notin\S$ in the disk of radius $R\ge
2$, we have an inequality
\begin{equation}\label{lower-bound}
    |p(t)|\ge 2^{-O(d)}(r/R)^d,\qquad r=\dist(t,\S),\ |t|\le R.
\end{equation}
\end{Lem}

\begin{proof}
For a polynomial $p=\sum_{j=0}^d p_jt^j$ of unit norm and degree
$d$, at least one of its $d+1$ coefficients must be at least
$1/(d+1)$ in the absolute value. Writing the Cauchy formula for the
corresponding derivative at the origin, we conclude that
\begin{equation*}
    p_j= \frac1{2\pi\iu}\oint_{|t|=1}\zeta^{-j}\,p(\zeta){\d
    \zeta}{},\qquad\text{hence}\quad
    \frac1{2\pi}\int_{|t|=1}|p(\zeta)|\,\d s\ge \frac1{d+1}.
\end{equation*}
By the mean value theorem there exists a point $t_*$ on the unit
circle, such that $|p(t_*)|\ge\frac1{d+1}$.

Let $K\subset\C$ be the disk of radius $2R$, centered at the origin.
The polynomial $p$ can be written as the product,
\begin{equation}\label{double-prod}
    p(t)=\alpha\prod_{|t_j|>2R}(t-t_j)\prod_{|t_j|\le 2R}\bigl(1-\tfrac t{t_j}\bigr),\qquad
    \alpha\in\C,
\end{equation}
where $t_j$ are roots of the polynomial, counted with their
multiplicity. Denote by $d_\pm$ the number of roots of $\S$ inside
(resp., outside) the disk $\{|t|\le 2R\}$, counted with
multiplicities so that $d_-+d_+=d$.

A lower bound for $|\alpha|$ follows from the inequality
$|p(t_*)|=1/(d+1)$. For $t=t_*$ the first product in
\eqref{double-prod} is \emph{majorized} by $(2R+1)^{d_-}$, while the
second product does not exceed $(4/3)^{d_+}$, since $|t|/|t_j|\le
1/(2R-1)\le 1/3$ for all terms in it. Altogether we conclude that
\begin{equation*}
    |\alpha|\ge\frac1{d+1}\cdot\frac{(3/4)^{d_+}}{(2R+1)^{d_-}}.
\end{equation*}

This implies the lower bound for $|p(t)|$ with $|t|\le R$:
\begin{equation*}
    |p(t)|\ge |\alpha|\,r^{d_-}(1/2)^{d_+}.
\end{equation*}
Combining these two estimates, we conclude that
\begin{equation*}
    |p(t)|\ge \frac1{d+1}\cdot
    \biggl(\frac{r}{2R+1}\biggr)^{d_-}(3/8)^{d_+}\ge 2^{-O(d)}\biggl(\frac rR\biggr)^d
\end{equation*}
for all  $2\le R<+\infty$ and $R>r>0$.
\end{proof}

\begin{Rem}
One can give a different bound using the Cartan lemma
\cite{levin}. This lemma establishes an explicit lower bound for a
monic polynomial outside a union of disks of given total diameter.
For our purposes it is less convenient, since the polynomial in
question is not monic (and can actually be very far from it).
\end{Rem}

\subsection{Variation of argument of solutions along arcs}
For the next step, we give an upper bound for the variation of
argument of complex-valued solutions of a Fuchsian equation along a
circular arc  (or line segment) in terms of the slope of
the equation and the normalized length of the arc.

\begin{Lem}\label{lem:var-arg}
Let $D$ be a differential operator of order $k$ with coefficients of
degree $\le d$ and slope $S=\slope D$, and $\gamma$ a closed
circular arc or line segment disjoint with the singular locus $\S$,
which belongs to the disk of radius $R$ centered at the origin.

Then the variation of argument of any nonzero solution of the
homogeneous equation $Dy=0$ along the arc $\gamma$ is explicitly
bounded,
\begin{equation}\label{var-arg-bd}
    \Var\Arg y(t)|_\gamma \le k{S}L(R/r)^{O(d)}.
\end{equation}
where $L$ is the length of the arc $\gamma$ and $r=\dist(\gamma,\S)$.
\end{Lem}

\begin{proof}
Without loss of generality we may assume that the polynomial
coefficients of the operator $D$ are normalized as follows,
\begin{equation*}
    \|a_0\|=1,\qquad \|a_j\|\le {S},\quad j=1,\dots,k,
\end{equation*}
where ${S}=\slope D$.

On the disk of radius $R$ the absolute values of the non-principal
coefficients are bounded, $|a_j(t)|\le {S}R^d$. Restricting the
equation $Dy=0$ on the arc $\gamma$ parameterized by the arc-length,
we obtain a homogeneous linear ordinary differential equation. Dividing by $a_0|_\gamma$ and applying Lemma~\ref{lem:lower-bound} to this polynomial, this equation can be brought into the
monic form with explicitly bounded complex valued coefficients,
\begin{equation}\label{monic-eq}
\begin{gathered}
    y^{(k)}+A_1(s)\,y^{(k-1)}+\cdots+A_k(s)=0,\qquad s\in[0,L],\\
     |A_j(s)|\le 2^{O(d)}\cdot
    {S}R^{2d}/r^d\le {S}(R/r)^{O(d)},\ j=1,\dots,k.
\end{gathered}
\end{equation}
Applying the complex generalization of the de la Vall\'ee
Poussin nonoscillation theorem \cite[Corollary 2.7]{fields} to \eqref{monic-eq}, we
conclude that the variation of argument admits the bound
\begin{equation}\label{fields}
    \Var\Arg y(\cdot)|_\gamma\le O(1)\cdot k
    L A,\qquad A=\max_{j=1,\dots,k}\max_{s\in\gamma}|A_j(s)|.
\end{equation}
Substituting the bounds \eqref{monic-eq} into this inequality, we
conclude with the bound \eqref{var-arg-bd}.
\end{proof}

The assertion of Lemma~\ref{lem:var-arg} allows one to place an explicit
upper bound for the number of isolated zeros of any solution to a
linear equation with rational coefficients in a simply connected
polygonal domain away from the singular locus. Indeed, one can
estimate the variation of argument of any solution along each arc
bounding the polygonal domain, and apply the argument principle.

In addition to the explicit bound on the variation of argument along arcs
distant from the singular points, one can slightly modify the argument
above to give a bound on the variation of argument along a small circle
centered at a \emph{Fuchsian singular point} as well.

\begin{Lem}
Let $t_0\in\S$ be a Fuchsian singular point of a linear operator
with polynomial coefficients of order $k$, degree $d$ and slope
${S}$, which belongs to the disk of radius $R\ge 2$ and is at least
$r$-distant from all other points $t_1,\dots,t_d\in\S$.

Then the variation of argument of any solution $y(t)$ along the
circle $\gamma_\rho=\{|t-t_0|=\rho\}$ is bounded,
\begin{equation}\label{eq:var-arg-small-circle}
    \Var\Arg y(\cdot)|_{\gamma_\rho}\le k{S}(R/r)^{O(d)}
\end{equation}
uniformly over all sufficiently small $\rho>0$.
\end{Lem}

\begin{proof}
Note that the parallel translation $p(\cdot)\mapsto
(\operatorname{Tr}_Rp)(\cdot)=p(\,\cdot\,+R)$ is a linear operator
of explicitly bounded norm on the space of polynomials of degree
$\le d$. Indeed, by the Taylor formula,
\begin{equation*}
    \|\operatorname{Tr}_Rp\|\le \sum_{i=0}^{d}\frac {R^i\, \|p^{(i)}\|}{i!}\le
    \|p\|\cdot \sum_0^d \frac{(Rd)^i}{i!}\le2^{O(d)}R^d\cdot\|p\|.
\end{equation*}
Thus without loss of generality we may assume that the Fuchsian
singularity $t_0=0$ is at the origin, and the slope of the new
differential equation is at most $2^{O(d)}{S}R^d$, where ${S}$ is
the original slope.

Since the point is Fuchsian, the equation can be written in the form
\begin{equation}\label{eq-fuchs}
    t^k q_0(t)\,y^{(k)}(t)+t^{k-1}q_1(t)\,y^{(k-1)}(t)+\cdots+q_k(t)y(t)=0,
\end{equation}
where $q_0,q_1,\dots,q_k$ are polynomials of degrees $\le
d-k,d-k+1,\dots,d$ respectively, see \cite[Proposition
19.19]{iy:lade}. Since the norm is multiplicative, and $\|t^j\|=1$
we conclude that $\|q_j\|/\|q_0\|\le 2^{O(d)}{S}R^d$ (see above). By
Lemma~\ref{lem:lower-bound}, $|q_0(t)|\ge 2^{-O(d)}r^d$ and,
obviously, $|q_j(t)|\le\|q_j\|$ for all sufficiently small
$|t|=\rho$.

Letting $t=\rho s$ and noting that $\frac{\d^k}{\d
t^k}=\rho^{-k}\frac{\d^k}{\d s^k}$, we rewrite the equation above in
the form
\begin{equation}\label{eq-fuchs2}
    s^k q_0(\rho s)\,y^{(k)}(s)+s^{k-1}q_1(\rho s)\,y^{(k-1)}(s)+\cdots+q_k(\rho s)y(s)=0,
\end{equation}
with $|s|=1$. Notice that the rescaling  of $\frac{\d^k}{\d t^k}$
precisely cancels out with the rescaling of the vanishing terms
$t^k$ in the coefficients. This circumstance is unique to Fuchsian
differential operators.

After dividing the equation above by the leading term (nonzero
and bounded from below), we obtain a monic linear ordinary
differential equation on the unit circle $|s|=1$ with
coefficients not exceeding an explicit constant
$A=2^{O(d)}{S}(R/r)^d$. The bound for the variation of argument
along the circle again follows from \eqref{fields}.
\end{proof}

\subsection{Counting zeros in topological annuli}
In this section we recall an ``\emph{argument principle}'' for
multivalued functions. This principle asserts that the number of
zeros of a solution of a Fuchsian differential operator in an
annulus is explicitly bounded, provided that the monodromy of the
operator has all eigenvalues on the unit circle and the annulus is
at a positive distance from the singular locus. The proof is based
on an idea due to Petrov (the so called ``Petrov's trick'').  The
result we need appears (in slightly less general form) in \cite{mrl-96}.

Let $K$ be an \emph{annulus} bounded by two disjoint circles
$C_{1,2}$ with centers on the real axis, and $D$ a monic linear
operator with coefficients holomorphic in the closure $\overline K$
and real on the real axis $K\cap\R$ (we will refer to such operators
as \emph{real}). In this case Lemma~\ref{lem:var-arg} asserts that
there exists an upper bound $B=B(D,K)$ for the variation of argument
of any solution of the homogeneous equation $Dy=0$ in $K$ along each
of the boundary arcs of $\partial K$,
\begin{equation}\label{var-arg-ann}
    \Var\Arg y(\cdot)|_{C_{1,2}}\le B=B(D,K).
\end{equation}

\begin{Lem}\label{lem:annulus-count}
If the monodromy of a \emph{real} differential operator $D$ along
the equator of an annulus $K$ symmetric around the real axis has all eigenvalues on the
unit circle, then the number of zeros of any solution in $K$ is
explicitly bounded,
\begin{equation}\label{annulus-bd}
    \cN(D,K)\le (2k+1)(2B+1),
\end{equation}
where $k$ is the order of the operator and $B=B(D,K)$ the bound from
\eqref{var-arg-ann}.
\end{Lem}

\begin{proof}
This Lemma almost literally coincides with Theorem~2 from
\cite{mrl-96}. To simplify the exposition, we first make a real
conformal automorphism of the complex line $\C P^1$, which
transforms the given annulus $K$ into an annulus $K'$ bounded by two
\emph{concentric} circles centered at the origin, and then pass to
the logarithmic chart $\zeta$ on the universal covering over $K'$.
As a result, we arrive at the following problem.

Let $\^\Pi$ be the ``vertical'' strip $c_-\le\Re\zeta\le c_+$ and
$D$ a linear ordinary operator with coefficients holomorphic in
$\^\Pi$, $2\pi\iu$-periodic in this strip and real on $\^\Pi\cap\R$.
Let $B$ be the upper bound for the variation of argument of any
solution of the equation $Dy=0$ along the segment of length $2\pi$
on any of the lines $\Re\zeta=c_\pm$. Because of the periodicity,
the shift operator $\mathbf M:y(\cdot)\mapsto y(\,\cdot\,+2\pi \iu)$
is an automorphism of the space of solutions of the equation,
corresponding to the monodromy of the original equation in the
annulus.

We need to show that if all eigenvalues of $\mathbf M$ are of unit
modulus, $|\mu|=1$, then the number of zeros of any solution in
the rectangle $\Pi=\^\Pi\cap\{|\Im \zeta|\le 2\pi\}$ is at most
$(k+1)(2B+1)$; this is equivalent to the assertion of the Lemma.

By the argument principle, it is sufficient to place an upper bound
on the variation of argument of any solution $f$ of the equation
$Df=0$ along the perimeter of the rectangle, which consists of two
vertical segments on the lines $\Re\zeta=c_\pm$ and two horizontal
segments $$S_\pm=\{\Im\zeta=\pm2\pi,\ c_-\le\Re\zeta\le c_+\}.$$ The
contribution of the vertical segments is explicitly given by the
assumptions of the lemma; it does not exceed $2B$. It remains to
estimate the variation of arguments along the horizontal segments.
Each increment of $\Arg f(\zeta)$ by an angle $\pi$ or more on a
connected interval of $\zeta$ implies, by an intermediate value
theorem, that at some point $\zeta_0$ on this interval the imaginary
part $\Im \bigl(\mu f(\zeta_0)\bigr)$ vanishes, for an arbitrary
choice of the complex number $\mu\ne0$. This observation implies the
inequality
\begin{equation*}
    \Var\Arg f|_{S_\pm}\le \pi(\nu_\pm +1),
    \qquad\nu_\pm=\{\zeta\in S_\pm\:\Im\bigl(\mu^{\mp1} f(\zeta)\bigr)=0\}
\end{equation*}
valid for any complex number $\mu\ne0$. In general, one has no explicit control over the imaginary part of a solution along the segments
$S_\pm$. However, if $f$ is real on the real segment
$S_0=\Pi\cap\R$, then by the Schwarz symmetry principle,
$f(\zeta+2\pi \iu)=\overline{f(\zeta-2\pi\iu)}$ for $\zeta\in\S_0$.
Therefore for any $\mu$ with $|\mu|=1$,
\begin{equation*}
    \Im (\mu^{\mp1}f)|_{S_\pm}=\Im (\mu^{\mp1}\mathbf
    M^{\pm1}f)|_{S_0}=\mathbf P_\mu f|_{S_0},
\end{equation*}
where $\mathbf P_\mu$ is the \emph{Petrov operator}
\begin{equation}\label{petrov}
    \mathbf P_\mu=\frac1{2\iu}\, (\mu^{-1}\mathbf M-\mu\mathbf M^{-1}).
\end{equation}
This yields an inequality relating the number $N(f,\Pi)$ of
\emph{complex} zeros of a real solution $f$ in the rectangle $\Pi$
to the number of \emph{real} zeros of the function $f'=\mathbf
P_\mu f$ on $S_0\subset\Pi$: for any $\mu$ with $|\mu|=1$,
\begin{equation}\label{petrov-rolle}
    N(f,\Pi)\le \frac1{2\pi }\biggl(2B +2\pi \bigl(N(f',S_0)\bigr)+1\biggl)\le
    (2B+1)+N(f',\Pi).
\end{equation}
This inequality is a complex analog for the Rolle inequality
relating the number of zeros of a real smooth function and its
derivative. The inequality remains valid also in the case when
$f'=\mathbf P_\mu f\equiv0$, if we set $N(0,\Pi)=0$. This follows
from the symmetry of the function $f$ (the contributions of
$\Var\Arg f$ along the segments $S_\pm$ cancel each other).

To complete the proof of the Lemma, we show that the product
(composition) of operators $\mathbf P_D=\prod_i\mathbf
P_{\mu_i}^{m_i}$ (where $\mu_i$ range over the  eigenvalues of the
monodromy operator $\mathbf M$ and $m_i$ over the corresponding
multiplicities) vanishes on the linear space of solutions of the
operator $Dy=0$ in $\^\Pi$. This assertion is an analog of the
Hamilton--Cayley theorem for the Petrov operators.

Indeed, suppose $\mathbf M$ is a Jordan block of size $m$ with an
eigenvalue $\mu$ and let $0\subset V_1\subset\cdots\subset
V_{\mu_i}\subset0$ denote the corresponding flag. Then $\mathbf M$
acts on $V_{j+1}/V_j$ as multiplication by $\mu$, and $\mathbf
P_{\mu}$ acts by zero. Thus $\mathbf P_{\mu}V_{j+1}\subset V_j$ and
it follows that $\mathbf P_{\mu}^{m}=0$. It remains to note that
Petrov operators with different numbers $\mu$  commute with each other.

Now the proof of the Lemma for functions real on $\R$ is
straightforward: if $Df=0$, then $\mathbf P_Df=0$, and by
\eqref{petrov-rolle}, the number of isolated zeros of $f$ does not
exceed $k(2B+1)$, where $k$ is the order of $D$, equal to the number
of the Petrov operators forming the product $\mathbf P$.

Finally consider the case where $f$ is a general complex valued
solution.  By the same argument as above, the total number of zeros
of $f$ in $\Pi\cap\{|\Im \zeta|\le 2\pi\}$ does not exceed
\begin{equation} \label{eq:var-non-real}
   B+\tfrac{1}{2\pi}\(\Var\Arg f|_{S_+}+\Var\Arg f|_{S_-}\).
\end{equation}
Assume that $\Im f|_{S_\pm}\not\equiv 0$ (if this assumption fails,
consider $\Re f|_{S_\pm}$ instead). Let $f_\pm$ denote the function
which agrees with $\Im (\mathbf M^\pm f)$ on $S_0$ and is extended
analytically onto $\Pi$. The values of $f^\pm$ on $S_0$ correspond
to the imaginary part of $f$ on $S^\pm$ respectively.

Since $D$ is a real operator, $f^\pm$ is again a solution of the
equation $Dy=0$.  These functions are real on $S_0$, the previous
considerations apply, and we conclude that
\begin{equation*}
   \Var\Arg f|_{S_\pm}\le 2\pi N(f^\pm,S_0)\le k(2B+1).
\end{equation*}
Combining this with \eqref{eq:var-non-real}
we finally obtain the bound
\eqref{annulus-bd}.
\end{proof}

\subsection{Symmetrization}
In order to apply Lemma~\ref{lem:annulus-count}, to a given annulus,
one needs the assumption that the coefficients of the operator $D$
are real on the real axis, eventually after the rotation of the
axis. In this section we show that this additional assumption can in
fact always be achieved, if the operator $D$ is obtained from a
system of linear differential equations \eqref{eq:ls}. This
construction will be referred to as the \emph{symmetrization} of the
operator relative to a given axis.

Let $f\in\mathscr O(U)$ be a function holomorphic in an open domain
$U\subset\C$. The function $f^\dag$, defined in the domain $U^\dag$
by the formula
\begin{equation}\label{symm-f}
    f^\dag(t)=\overline{f(\bar t)},\qquad t\in U^\dag=\{\bar t\:
    t\in U\},
\end{equation}
(the bar stands for the complex conjugacy) will be called a
\emph{reflection} of $f$ in the real axis. Obviously, $f^\dag$ is
holomorphic, $f^\dag\in\mathscr O(U^\dag)$. For polynomials the
reflection consists replacing all (constant) coefficients by their
complex conjugates:
\begin{equation*}
    p(t)=\sum_{k=0}^m c_k\,t^k \iff p^\dag(t)=\sum_{k=0}^m \bar
    c_k\,t^k,\qquad c_0,\dots,c_m\in\C,
\end{equation*}
and it extends to the rational functions in the obvious way.

Let $U=U^\dag\subseteq\C$ be a domain symmetric with respect to the
real axis $\R\subset\C$. A $\C$-linear subspace $L\subseteq\mathscr
O(U)$ of functions defined in $U$, is called \emph{symmetric} (with
respect to the real line), if together with each function
$f\in\mathscr O(U)$ it contains the function $f^\dag$. A symmetric
subspace is closed under taking the real/imaginary part on $U\cap\R$:
for any $f\in L$ the functions $\Re f|_{U\cap\R}$ and $\Im f_{U\cap
\R}$ can be analytically continued from the subset of real points of
$U$ onto the entire domain $U$ using the formulas $\Re
f=\tfrac12(f(t)+f^\dag(t))$, $\Im f=\tfrac1{2\iu}(f(t)+f^\dag(t))$
respectively. A symmetric linear space (in case it is
finite-dimen\-sion\-al) always admits a basis of functions real on
$U\cap\R$. A monic linear ordinary differential operator
$D\in\mathscr O(U)[\partial]$ defines a symmetric space of
solutions, if and only if its coefficients are real on $\R$, i.e.,
$D\in\mathscr O_\R(U\cap\R)[\partial]$. Such operators will be
referred to as \emph{real} (on $\R$).

Any finite-dimensional subspace $L\subset\mathscr O(U)$ of functions
can be embedded into the symmetric completion $L^\ominus=L+L^\dag$
by simply adjoining the reflections $f^\dag$ of all elements $f\in
L$. As a consequence, each linear operator $D\in\mathscr
O(U)[\partial]$ can be considered as a ``right factor'' of a real
operator $D^\ominus\in\mathscr O_\R(U\cap\R)[\partial]$ of higher
order with real coefficients, which vanishes on all functions from
$L^\ominus$. However, the slope $\slope D^\ominus$ of the operator
may increase uncontrollably (relatively to $\slope D$) after such
symmetrization.

\begin{Ex}
The function $f(t)=\iu +\frac1{100}t^{100}$ satisfies a first order
equation
\begin{equation*}
    y'-\frac{t^{99}}{\iu+\frac1{100}t^{100}}\,y=0
\end{equation*}
which has very small (in the absolute value) coefficient in the disk
$|t|<\tfrac12$ and the slope $\slope D$ very close to $1$. On the
other hand, the symmetrization $L^\ominus$ is generated by the
functions $1$ and $t^{100}/100$ and is defined by the second order
equation
\begin{equation*}
    y''-\frac{99}t \,y'=0.
\end{equation*}
This equation has uniformly large coefficients in the same disk (and
a pole at the origin). Even the slope of the corresponding second
order operator is relatively large, $\slope D^\ominus=99$.
\end{Ex}

However, at the level of linear first order \emph{systems} the
symmetrization can be achieved very easily. Let
$\Omega^\dag\in\clF$ be the rational matrix 1-form built from
the 1-form $\Omega\in\clF$ on $\C$ by the reflection in the
real axis as above (replacing all coefficients of the entries
$\omega_{ij}$ by their conjugates). Then each solution of the
reflected system $\d X^\dag=\Omega^\dag\cdot X^\dag$ can be
obtained by reflection of the corresponding solution of the
initial system $\d x=\Omega x$ \eqref{eq:ls} and vice versa.
Therefore the symmetric linear space $L^\ominus$ spanned by all
components of a solution of the system
\begin{equation}\label{symmetrization}
    \d y=\Omega^\ominus y,\qquad
    \Omega^\ominus=\begin{pmatrix}\Omega&\\&\Omega^\dag\end{pmatrix},\ y\in\C^{2n},
\end{equation}
includes among its solutions the symmetrization of the solution of
the initial system. Reducing the doubled system
\eqref{symmetrization} to a scalar equation, one obtains an explicit
expression for the operator $D^\ominus$.

Clearly, one can replace the real axis $\R\subseteq\C$ by any (real
one-dimensional) line $\ell$ in the complex plane $\C$. Any such
line can be mapped to the real axis by a suitable affine
transformation $\f\:\C\to\C$ (in fact, rotation and translation
suffice). The corresponding reflection operator replacing
\eqref{symm-f}, can be explicitly written as follows,
\begin{equation*}
    f\mapsto f^\dag_\ell=(f\circ\f)^\dag\circ \f^{-1}.
\end{equation*}

Let now $U$ be an arbitrary concentric annulus free from singular
points of a Fuchsian system $\Omega\in\clF$. Among all systems
obtained by reflection of the matrix 1-form $\Omega$ in different
axes $\ell$ passing through the center, one can always find a
system (denote it again by $\Omega^\dag$) whose singular points
(except for the singular point at the center) will be sufficiently
far away from the singularities of $\Omega$. More precisely, we can
always guarantee that the matrix form $\Omega^\ominus$ of the
symmetrized system \eqref{symmetrization} satisfy the inequality
\begin{equation}\label{double-carpet}
    R_\flat(\Omega^\ominus)\le R_\flat(\Omega)^{\nu},\qquad
    \nu=O(n+m),
\end{equation}
with an explicit constant in $O(\cdot)$. Indeed, the norms of the
residues of the symmetrized system are majorized by the norms of
residues of $\Omega$ (residues of the matrix form $\Omega^\dag$ have
the same norm as those of $\Omega$). As for the distances between
the singular locus $\S=\{\tau_j\}$ and its mirror image
$\S^\dag=\S^\dag_\ell$, we can show that that they can be offset by
a suitable rotation of the axis of the symmetry $\ell\subset\C$.
Indeed, rotation of $\ell$ by an angle $\alpha$ results in the
rotation of $\S^\dag_\ell$ by the angle $2\alpha$. Yet for any two
$m$-point sets $A,B\subset\{\tau=1\}$ on the unit circle, one can
always find a rotation $\rho$ of the circle so that the arclength
distance between the points of $A$ and $\rho(B)$ will not exceed
$2\pi/m^2$ (a simple consequence of the pigeonhole principle). This
shows that even in the worst case when the singular locus belongs to
a single circle, the distance between $\S$ and $\S^\dag$ can always
be assumed bounded from below. The general case of arbitrary
singular locus can be immediately reduced to this case by
controlling the angular distance between the loci and their distance
from the center.

\subsection{Demonstration of Theorem~\ref{thm:main}}
Consider an arbitrary Fuchsian linear system \eqref{eq:ls} from the
Fuchsian class $\clF$ with the matrix form $\Omega$ the
corresponding system of linear equations in the affine chart $t$ as
in \eqref{eq:lst}.

This system after reduction to a linear differential equation of
order $n$ takes the form \eqref{eq:linear} with the coefficients
$a_0,\dots,a_n\in\C[t]$ of the corresponding differential operator
$D$ of degree $d\le n^2m$ and the slope not exceeding
$R^{2^{O(n^2m)}}$, $R=R_\flat$, by Lemma~\ref{lem:scalar-eq}.
Solutions of this equation are analytic outside the singular locus
$\S=\{\tau_0,\dots,\tau_m\}$, though the leading coefficient $a_0$
vanishes at some other points as well (apparent singularities of the
equation).

\begin{figure}\label{fig:slits}
\includegraphics[width=0.5\hsize, bb = 0 0 489 489]{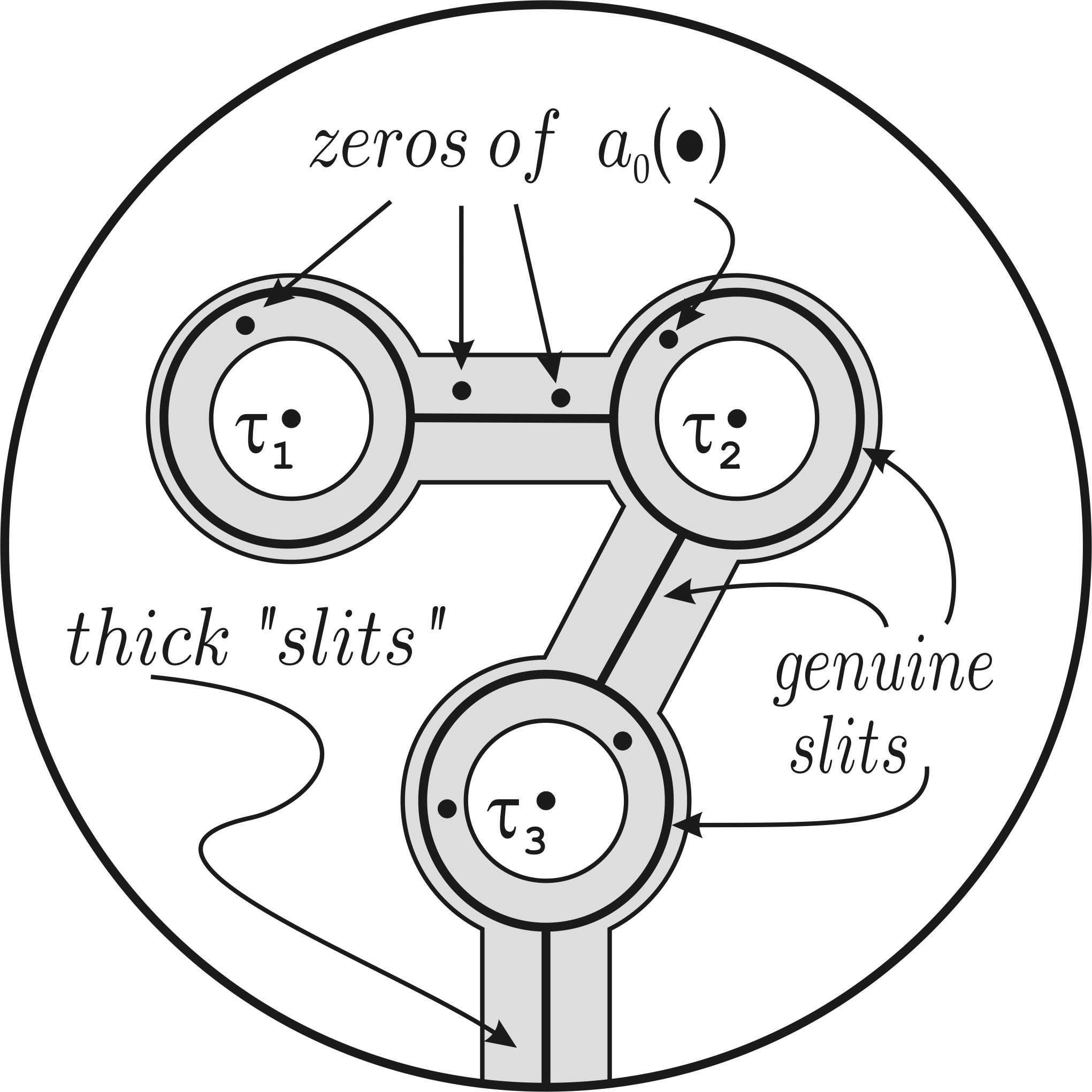}
\caption{``Thick slits'' and genuine slits.}
\end{figure}

Moreover, for any singular point $\tau\in\S$, one can find a real
line $\tau\in\ell\subset\C$ (``axis of symmetry'') passing through
this point, such that the linear operator $D^\ominus$ obtained by
symmetrization of $D$ in this axis has the order $2n$ and the slope
bounded by the same double exponential expression $R^{2^{O(n^2m)}}$.
Indeed, one has to apply Lemma~\ref{lem:scalar-eq} to the
symmetrized system \eqref{symmetrization}, in which $\Omega^\dag$ is
the result of reflection of $\Omega$ in $\ell$ and note that the
carpeting function will be replaced by a polynomially equivalent one
\eqref{double-carpet}. Abusing the language, we will say that the
operator $D$ admits an axis of symmetry.

Consider the system of ``thick slits'' of the disk
$K_R=\{|t|\le\R\}$ as follows (see Fig.~\ref{fig:slits}):
\begin{enumerate}
 \item annuli $\{\tfrac14\le |t-\tau_j|\le \tfrac12\}$ around each singular
 point $\tau_j\in\S$;
 \item rectangular strips of width $\tfrac12$ along all segments
 $[\tau_i,\tau_j]$ (these segments should be axes of symmetry for these rectangles);
 \item rectangular strips of with $\tfrac12$ along the
 segment connecting $\tau_j$ with the exterior circle $\partial K_R=\{|t|=R\}$.
\end{enumerate}
By elementary geometric considerations it is clear that one can
construct a system of genuine slits of $K_R\ssm\S$ as follows:
\begin{enumerate}
 \item these slits are only circular arcs and line segments;
 \item they subdivide $K_R\ssm\S$ into simply connected polygonal
 domains and punctured disks;
 \item each slit (circular or straight) is at least $r$-distant from
 the null locus of the leading coefficient $V=\{a_0=0\}$, where $1/r=O(n^2m)$;
 \item the total length of the slits is bounded by $O(m^2R)$.
\end{enumerate}
Indeed, even if all $d=m^2n$ points from $V$ belong to a given
annulus or a rectangular strip of width $O(1)$, one can always find
a circle (resp., line segment) which is at least $O(1/d)$-distant
from $V$.

We can now assemble all bounds given by Lemma~\ref{lem:scalar-eq},
Lemmas~\ref{lem:var-arg} in common terms relevant to the Fuchsian
class $\clF$ and the carpeting function $R=R_\flat$ on it. The
variation of arguments of any linear combination of solutions of the
system is explicitly bounded along all slits constructed above, and
using the argument principle for simply connected pieces (resp.,
Lemma~\ref{lem:annulus-count} for annuli), we arrive at the explicit
bounds for the counting function.

Obviously, in all products \eqref{var-arg-bd},
\eqref{eq:var-arg-small-circle} and hence in \eqref{annulus-bd}, the
double exponential term ${S}=R^{2^{O(d)}}$ absorbs all other terms
without changing the asymptotic behavior of the entire products,
which ultimately proves the bound asserted in
Theorem~\ref{thm:main}. \qed

\bibliographystyle{cdraifplain}
%\def\MR{}
%\bibliography{tangent16,2003,2007}
\def\BbbR{$\mathbb R$}\def\BbbC{$\mathbb
  C$}\providecommand\cprime{$'$}\providecommand\mhy{--}\font\cyr=wncyr8\def\Bb%
bR{$\mathbb R$}\def\BbbC{$\mathbb
  C$}\providecommand\cprime{$'$}\providecommand\mhy{--}\font\cyr=wncyr9\def\Bb%
bR{$\mathbb R$}\def\BbbC{$\mathbb
  C$}\providecommand\cprime{$'$}\providecommand\mhy{--}\font\cyr=wncyr9\def\cp%
rime{$'$}
\def\bysame{\leavevmode ---------\thinspace}
\makeatletter\if@francais\providecommand{\og}{<<~}\providecommand{\fg}{~>>}
\else\gdef\og{``}\gdef\fg{''}\fi\makeatother
\def\cdrandname{\&}
\providecommand\cdrnumero{no.~} \providecommand{\cdredsname}{eds.}
\providecommand{\cdredname}{ed.}
\providecommand{\cdrchapname}{chap.}
\providecommand{\cdrmastersthesisname}{Memoir}
\providecommand{\cdrphdthesisname}{PhD Thesis}

\end{document}